\documentclass[3p]{elsarticle}
\usepackage{amsmath}
\usepackage{amsfonts,amssymb}

\newtheorem{theorem}{Theorem}{\bf}{\it}
\newtheorem{prop}{Proposition}[section]
{\bf}{\it}
{\bf}{\rm}
\newtheorem{lemma}{Lemma}{\bf}{\it}
{\bf}{\it}
{\bf}{\it}
{\bf}{\it}
{\bf}{\rm}{\rm}

\newcommand{\zr}{\ltimes}
\newcommand{\Real}{\mathbb{R}}
\newcommand{\Co}{\mathbb{C}}

\newcommand{\g}{\mathfrak{g}}
\newcommand{\h}{\mathfrak{h}}

\newcommand{\so}{\mathfrak{so}}

\newcommand{\simil}{\mathfrak{sim}}

\def\u{\mathfrak{u}}
\def\sp{\mathfrak{sp}}
\def\spp{\mathfrak{sp}}

\newcommand{\z}{\mathfrak{z}}
\newcommand{\R}{\mathcal{R}}
\newcommand{\M}{\mathcal{M}}

\def\P{\mathcal{P}}
\newcommand{\id}{\mathop\text{\rm id}\nolimits}
\newcommand{\spa}{\mathop\text{{\rm span}}\nolimits}
\newcommand{\Hom}{\mathop\text{\rm Hom}\nolimits}
\newcommand{\End}{\mathop\text{\rm End}\nolimits}

\newcommand{\ad}{\mathop\text{\rm ad}\nolimits}

\newcommand{\pr}{\mathop\text{\rm pr}\nolimits}
\newcommand{\Op}{\mathop\text{\rm Op}\nolimits}
\newcommand{\Sim}{\mathop\text{\rm Sim}\nolimits}

\def\g{\mathfrak{g}}
\def\h{\mathfrak{h}}
\def\so{\mathfrak{so}}

\def\sp{\mathfrak{sp}}

\def\u{\mathfrak{u}}

\def\spa{\mathop\text{{\rm span}}\nolimits}
\def\Hom{\mathop\text{\rm Hom}\nolimits}
\def\pr{\mathop\text{\rm pr}\nolimits}

\def\End{\mathop\text{\rm End}\nolimits}
\def\id{\mathop\text{\rm id}\nolimits}

\def\ad{\mathop\text{\rm ad}\nolimits}

\def\Im{\mathop\text{\rm Im}\nolimits}

\def\Re{\mathop\text{\rm Re}\nolimits}

\def\SO{\text{\rm SO}}

\def\M{{\cal M}}

\def\zr{\ltimes}

\def\spa{\text{\rm span}}
\def\Hom{\text{\rm Hom}}

\def\End{\text{\rm End}}

\def\Real{\mathbb{R}}
\def\Co{\mathbb{C}}
\def\g{\mathfrak{g}}
\def\h{\mathfrak{h}}

\def\so{\mathfrak{so}}

\def\u{\mathfrak{u}}

\def\sp{\mathfrak{sp}}

\def\z{\mathfrak{z}}
\def\R{\mathcal{R}}

\def\H{\mathbb {H}}

\def\id{\text{\rm id}}

\def\pr{\text{\rm pr}}

\def\im{\text{\rm Im}}
\def\re{\text{\rm Re}}

\def\ad{\text{\rm ad}}

\def\Mat{\mathop\text{\rm Mat}\nolimits}

\begin{document}

\begin{frontmatter}

\title{Holonomy algebras of pseudo-hyper-K\"ahlerian manifolds of index 4}

\author{Natalia I. Bezvitnaya}
\address{Department of Mathematics and Statistics, Faculty of Science, Masaryk University in Brno,
Kotl\'a\v rsk\'a~2, 611~37 Brno, Czech Republic\\
\ead{bezvitnaya@math.muni.cz}}

\begin{abstract} The holonomy algebra of a pseudo-hyper-K\"ahlerian manifold of signature
$(4,4n+4)$ is a subalgebra of $\sp(1,n+1)$. Possible holonomy
algebras of these manifolds are classified.
 Using this, a new
proof of the classification of simply connected
pseudo-hyper-K\"ahlerian symmetric spaces of index 4 is obtained.
\end{abstract}


\end{frontmatter}

\section{Introduction}
The classification of  holonomy algebras of Riemannian manifolds
is well known and it has a lot of applications both in geometry
and physics, see e.g. \cite{Ber,Besse,Bryant,Gibbons09,Jo}. Lately
the theory of pseudo-Riemannian geometries has been steadily
developing. In particular, a classification of holonomy algebras
of pseudo-Riemannian manifolds is an actual problem of
differential geometry. It is solved only in some cases. The
difficulty appears if the holonomy algebra preserves a degenerate
subspace of the tangent space. Classification of holonomy algebras
of Lorentzian manifolds  is obtained in
\cite{Schell,BB-I,Leistner,Gal5,ESI}; classification of holonomy
algebras   of pseudo-K\"ahlerian manifolds of index 2 is achieved
in  \cite{GalDB}.
 These algebras are contained in  $\so(1,n+1)$ and $\u(1,n+1)\subset\so(2,2n+2)$, respectively.
There are partial results for holonomy algebras of
pseudo-Riemannian manifolds of signature $(2,n)$ and $(n,n)$
\cite{I99,BI97,GT01}.  More details can be found in the recent
review \cite{IRMA}.

In \cite{Quat} holonomy algebras of
pseudo-quaternionic-K\"ahlerian manifolds with non-zero scalar
curvature are classified. These algebras $\g$ are contained in
$\sp(1)\oplus\sp(r,s)$ and they contain $\sp(1)$. If $s\neq r$,
then $\g$  is irreducible. If $s=r$, then $\g$ may preserve a
degenerate subspace of the tangent space, in this case there are
only two possibilities for $\g$. This strong result follows mainly
from the inclusion $\sp(1)\subset\g$.

Recall that a pseudo-hyper-K\"ahlerian manifold  is a
pseudo-Riemannian manifold $(M,h)$ together with three parallel
$g$-orthogonal complex structures $I_{1},I_{2},I_{3}$ that satisfy
the relations $I_{1}^{2}=I_{2}^{2}=I_{3}^{2}=-\id$,
$I_{3}=I_{1}I_{2}=-I_{2}I_{1}$. Any such manifold has signature
$(4r,4s)$, $r+s>1$, and its holonomy algebra $\g$  is contained in
$\sp(r,s)$. Conversely, any simply connected pseudo-Riemannian
manifold with such holonomy algebra is pseudo-hyper-K\"ahlerian.
Note that any  pseudo-hyper-K\"ahlerian manifold
 is also  pseudo-quaternionic-K\"ahlerian and it has zero scalar curvature.

In the present paper we classify all possible holonomy algebras
$\g\subset\sp(1,n+1)$ of pseudo-hyper-K\"ahlerian manifolds of
signature $(4,4n+4)$, $n\geq 1$. For $n=0$ this classification is
obtained in \cite{44}. The main results is stated in Section
\ref{REZ}. Section \ref{D} is dedicated to their proofs. To prove
the classification theorem we use the fact that a holonomy algebra
$\g\subset\sp(1,n+1)$ is {a Berger algebra}, i.e. $\g$ is spanned
by the images of the algebraic curvature tensors $R\in\R(\g)$ of
type $\g$. Recall that $\R(\g)$ is the space of linear maps from
$\wedge^2\Real^{4,4n+4}$ to $\g$ satisfying the first Bianchi
identity. In \cite{B1} weakly irreducible subalgebras
$\g\subset\sp(1,n+1)$ containing a certain ideal $\Im\H$ (see
decomposition \eqref{spHpdec1}) are partially classified. Here we
find missing subalgebras, then we compute the spaces $\R(\g)$ for
each of these algebras and we check which $\g$ are Berger
algebras. Then we show that each weakly-irreducible Berger
subalgebra $\g\subset\sp(1,n+1)$ contains $\Im\H$. This gives the
classification of weakly-irreducible not irreducible Berger
subalgebras $\g\subset\sp(1,n+1)$. Remark that in this paper only
possible holonomy algebras are listed. We do not know if all these
algebras may appear as the holonomy algebras, to show this
examples of manifolds must be constructed. Since the most of the
previously known Berger algebras are realized as the holonomy
algebras, one may expect that the algebras obtained here can be
realized as the holonomy algebras of pseudo-hyper-K\"ahlerian
manifolds.

In \cite{AC01,KO07,KO08} simply connected pseudo-hyper-K\"ahlerian
symmetric spaces of index 4 are classified. Using the results of
this paper we give a new simple proof of this classification. We
use the known fact that a simply connected pseudo-Riemannian
symmetric space is uniquely defined (up to a homothety) by its
holonomy algebra $\g$ and an algebraic curvature tensor
$R\in\mathcal{R}(\g)$ such that the representation of $\g$ in
$\mathcal{R}(\g)$ annihilates $R$ and the image of $R$ spans $\g$.
In Section \ref{secSym} we describe all such pairs $(\g,R)$ for
$\g\subset\sp(1,n+1)$. We show that if a pseudo-hyper-K\"ahlerian
manifold $(M,h)$ of signature $(4,4n+4)$, $n\geq 1$ is locally
symmetric, then $n=2$ and we give explicitly the curvature tensor
and holonomy algebra of the obtained space. For the case of
signature $(4,4)$ the analogous result is obtained in \cite{44}.

{\it Acknowledgement.} I am grateful to D.~V.~Alekseevsky and Jan
Slov\'ak for useful discussions, help and support. I am thankful
to Bastian Brandes for pointing out a gap in the first version of
this paper and for draw my attention to the paper \cite{DR78}. The
author has been supported by the grants 201/09/H012 of the Grant
Agency of Czech Republic and CZ.1.07/2.3.00/20.0003 of the
Ministry of Education, Youth and Sports of Czech Republic.

\section{Preliminaries}\label{secPril}
Let $\mathbb{H}^{m}$ be an $m$-dimensional quaternionic vector
space. A pseudo-quaternionic-Hermitian metric $g$ on
$\mathbb{H}^{m}$ is a non-degenerate $\mathbb{R}$-bilinear map
$g:\mathbb{H}^{m}\times\mathbb{H}^{m}\to\mathbb{H}$ such that
$g(aX,Y)=ag(X,Y)$ and $\overline{g(Y,X)}=g(X,Y)$, where
$a\in\mathbb{H}$, $X,Y\in\mathbb{H}^{m}$. Hence,
$g(X,aY)=g(X,Y)\bar{a}$. There exists a basis $e_{1},...,e_{m}$ of
$\mathbb{H}^{m}$ and integers $(r,s)$ with $r+s=m$ such that
$g(e_{t},e_{l})=0$ if $t \neq l$, $g(e_{t},e_{t})=-1$ if $1 \leq t
\leq r$ and $g(e_{t},e_{t})=1$ if $r+1 \leq t \leq m$. The pair
$(r,s)$ is called  the signature of $g$. In this situation we
denote $\mathbb{H}^{m}$ by $\mathbb{H}^{r,s}$. The realification
of $\mathbb{H}^{m}$ gives us the vector space $\mathbb{R}^{4m}$
with the quaternionic structure $(i,j,k)$. Conversely, a
quaternionic structure on $\mathbb{R}^{4m}$, i.e. a triple
$(I_{1},I_{2},I_{3})$ of endomorphisms of $\mathbb{R}^{4m}$ such
that $I_{1}^{2}=I_{2}^{2}=I_{3}^{2}=-\id$ and
$I_{3}=I_{1}I_{2}=-I_{2}I_{1}$, allows us to consider
$\mathbb{R}^{4m}$ as $\mathbb{H}^{m}$. A
pseudo-quaternionic-Hermitian metric $g$ on $\mathbb{H}^{m}$ of
signature $(r,s)$ defines on $\mathbb{R}^{4m}$ the
$i,j,k$-invariant pseudo-Euclidian metric $\eta$ of signature
$(4r,4s)$, $\eta(X,Y)=\re g(X,Y)$, $X,Y\in\mathbb{R}^{4m}$.
Conversely, a $I_{1},I_{2},I_{3}$-invariant pseudo-Euclidian
metric on $\mathbb{R}^{4m}$ defines a
pseudo-quaternionic-Hermitian metric $g$ on $\mathbb{H}^{m}$,
$$g(X,Y)=\eta(X,Y)+i\eta(X,I_{1}Y)+j\eta(X,I_{2}Y)+k\eta(X,I_{3}Y).$$

We will identify $(1,i,j,k)$ with $(I_0,I_1,I_2,I_3)$,
respectively. The identification $\Real^{4r,4s}\simeq\H^{r,s}$
allows to multiply the vectors of $\Real^{4r,4s}$ by quaternionic
numbers.

The Lie algebra $\sp(r,s)$ is defined as follows
\begin{align*} \sp(r,s)&=\{f\in\so(4r,4s) |
[f,I_{1}]=[f,I_{2}]=[f,I_{3}]=0\}\\&=
\{f\in\End(\mathbb{H}^{r,s})|\, g(fX,Y)+g(X,fY)=0 \text{ for all }
X,Y\in\mathbb{H}^{r,s}\}.\end{align*} Denote by $\sp(1)$ the
subalgebra in $\so(4r,4s)$  generated by the $\mathbb{R}$-linear
maps $I_{1},I_{2},I_{3}$.

Clearly, the tangent space of a pseudo-hyper-K\"ahlerian manifold
$(M,h)$ at a point $x\in M$ one can identify with
$(\mathbb{R}^{4r,4s},\eta,I_{1},I_{2},I_{3})=(\mathbb{H}^{r,s},g)$.
Then the holonomy algebra of a  pseudo-hyper-K\"ahlerian manifold
is identified with a subalgebra $\g\subset\sp(r,s)$.

Let $(V,\eta)$ be a pseudo-Euclidean space and $\g\subset\so(V)$
be a subalgebra. The space of curvature tensors $\mathcal{R}(\g)$
of type $\g$ is defined as follows
$$\mathcal{R}(\g)=\{R\in\Hom(\wedge^{2} V,\g)\ |\
R(u, v)w+R(v, w)u+R(w, u)v=0 \ for \ all \
 u,v,w\in V\}.$$

 Denote by $L(\mathcal{R}(\g))$ the vector subspace of $\g$
 spanned by the elements $R(u, v)$ for all
 $R\in\mathcal{R}(\g)$ and $u,v\in V$.
{\it A subalgebra $\g\subset\so(r,s)$ is called a Berger algebra
if $L(\mathcal{R}(\g))=\g$.} From the Ambrose-Singer theorem it
follows that if $\g\subset\so(V)$ is the holonomy algebra of a
pseudo-Riemannian manifold, then $\g$ is a Berger algebra.
Therefore,  Berger algebras may be considered as the candidates to
the holonomy algebras.

Now we summarize some facts about quaternionic vector spaces. Let
$\mathbb{H}^{m}$ be an m-dimensional quaternionic vector space and
$e_{1},...,e_{m}$ a basis of $\mathbb{H}^{m}$. We identify an
element $X\in\mathbb{H}^{m}$ with the column $(X_{t})$ of the left
coordinates of $X$ with respect to this basis,
$X=\sum_{t=1}^{m}X_{t}e_{t}$. Let
$f:\mathbb{H}^{m}\to\mathbb{H}^{m}$ be an $\mathbb{H}$-linear map.
Define the matrix $\Mat_{f}$ of $f$ by the relation
$fe_{l}=\sum_{t=1}^{m}(\Mat_{f})_{tl}e_{t}$. Now if
$X\in\mathbb{H}^{m}$, then $fX=(X^{t}\Mat_{f}^{t})^{t}$ and
because of the non-commutativity of the quaternionic numbers this
is not the same as $\Mat_{f}X$. Conversely, to an $m\times m$
matrix $A$ of the quaternionic numbers we put in correspondence
the linear map $\Op{A}:\mathbb{H}^{m}\to\mathbb{H}^{m}$ such that
$\Op{A}\cdot X=(X^{t}A^{t})^{t}$. If
$f,g:\mathbb{H}^{m}\to\mathbb{H}^{m}$ are two $\mathbb{H}$-linear
maps, then $\Mat_{fg}=(\Mat_{g}^{t}\Mat_{f}^{t})^{t}$. Note that
the multiplications by the imaginary quaternionic numbers are not
$\mathbb{H}$-linear maps. Also, for $a,b\in\mathbb{H}$ holds
$\overline{ab}=\bar{b}\bar{a}$. Consequently, for two square
quaternionic matrices we have
$(\overline{AB})^{t}=\bar{B}^t\bar{A}^t$.

\section{Results}\label{REZ}

Let $(M,h)$ be a pseudo-hyper-K\"ahlerian manifold of signature
$(4,4n+4)$, $n\geq 1$. The tangent space to the manifold $(M,h)$
at a point $x\in M$ can be identified with the pseudo-Euclidean
space $(\Real^{4,4n+4},\eta,I_{1},I_{2},I_{3})$, where $\eta$ is a
pseudo-Euclidean metric on $\Real^{4,4n+4}$, $(I_{1},I_{2},I_{3})$
is the quaternionic structure on $\Real^{4,4n+4}$. This space one
can identify with the pseudo-quaternionic-Hermitian space
$(\H^{1,n+1},g)$, where $g$ is a pseudo-quaternionic-Hermitian
metric on $\H^{1,n+1}$.

The Wu Theorem \cite{Wu} allows to assume that the manifold
$(M,h)$ is locally indecomposable, i.e. locally it is not a
product of pseudo-Riemannian manifolds of positive dimensions.
This happens if and only if the holonomy algebra
$\g\subset\sp(1,n+1)$ of $(M,h)$ does not preserve any proper
non-degenerate subspace of $\Real^{4,4n+4}$. Such subalgebras
$\g\subset\sp(1,n+1)$ are called {\it weakly irreducible}. If the
holonomy algebra $\g\subset\sp(1,n+1)$ is irreducible, then
$\g=\sp(1,n+1)$ \cite{Ber,Bryant,Sch}. Thus we may assume that
$\g\subset\sp(1,n+1)$ is weakly irreducible and not irreducible.
In this case $\g$ preserves a four-dimensional isotropic
$I_1,I_2,I_3$-invariant subspace $W\subset \mathbb{H}^{1,n+1}$. We
fix a non-zero vector $p\in W$, then $W=\mathbb{H}p$. Let
$q\in\H^{1,n+1}$ be any isotropic vector such that $g(p,q)=1$.
Denote by $\H^n$ the $g$-orthogonal complement to $\H p\oplus\H q$
in $\H^{1,n+1}$. Let $e_1,...,e_n$ be a basis of $\H^n$ and let
$G$ be the corresponding Gram matrix of $g|_{\H^n}$, i.e.
$G_{ab}=g(e_a,e_b)$. Denote by $\sp(1,n+1)_{\mathbb{H}p}$ the
maximal subalgebra of $\sp(1,n+1)$ that preserves the quaternionic
isotropic line $\mathbb{H}p$, this Lie algebra has the matrix
form:
\begin{equation}\label{spHp}\sp(1,n+1)_{\H p}=
\left.\left\{\Op\left(\begin{array}{ccc}a&-(G\bar{X})^t&b\\0&\Mat_{A}&X\\0&0&-\bar
a\end{array}\right)\right|\begin{array}{c}a\in\H,\,A\in\sp(n),\\X\in\H^n,\,
\,b\in\Im \H\end{array}\right\}.\end{equation}

Here $\Op\M$ denotes the $\H$-linear endomorphism of $\H^{1,n+1}$
given by a matrix $\M$, see Section \ref{secPril}.

Notice that to find the matrix form of $\sp(1,n+1)_{\mathbb{H}p}$
acting in $\mathbb{R}^{4,4n+4}$ one can consider the basis
$p,I_{1}p,I_{2}p,I_{3}p,e_1,...,I_3e_n,q,I_{1}q,I_{2}q,I_{3}q$ of
$\mathbb{R}^{4,4n+4}$ and  it is enough to change each element
$c=c_0+c_1i+c_2j+c_3 k$ of the matrix from \eqref{spHp} to the
matrix $$ \left(\begin{array}{cccc} c_0&-c_1&-c_2&-c_3\\
c_1&c_0&c_3&-c_2\\
c_2&-c_3&c_0&c_1\\
c_3&c_2&-c_1&c_0\end{array}\right).$$

We denote the element from \eqref{spHp} by the quadruple
$(a,A,X,b)$. One can easily find the following Lie brackets:
\begin{align*}[(a,0,0,0),(a',0,X,b)]=&(a'a-aa',0,\bar{a}X,2\im ba),\\
[(0,0,X,0),(0,0,Y,0)]=& (0,0,0,2 \im g(X,Y)),\\
[(0,A,0,0),(0,B,X,0)]=&(0,[A,B]_{\sp(n)},A X,0),\end{align*} where
$a,a'\in \mathbb{H}$, $X,Y\in \mathbb{H}^{n}$, $A,B\in\sp(n)$,
$b\in \im\mathbb{H}$.

Thus we get the decomposition
\begin{equation}\label{spHpdec}\sp(1,n+1)_{\H
p}=\H\oplus\sp(n)\zr(\H^n\zr\Im\H).\end{equation} We may also
write $\H=\Real\oplus\spp(1)$, then
\begin{equation} \label{spHpdec1}\sp(1,n+1)_{\H
p}=\Real\oplus\sp(1)\oplus\sp(n)\zr(\H^n\zr\Im\H).\end{equation}

The isomorphism $\{(a,0,0,0)|\, a\in\Im\H\}\simeq\spp(1)$ is given
by $(a,0,0,0)\mapsto -a$.

The ideal $\H^n\zr\Im\H\subset\sp(1,n+1)_{\H p}$ is isomorphic to
the quaternionic Heisenberg Lie algebra. The Levi-Malcev
decomposition of $\sp(1,n+1)_{\H p}$ has the form
$$\sp(1,n+1)_{\H p}=\mathfrak{s}\zr\mathfrak{r},\quad
\mathfrak{s}=\spp(1)\oplus\sp(n),\quad
\mathfrak{r}=\Real\zr(\H^n\zr\Im\H),$$ where $\mathfrak{s}$ is a
semisimple subalgebra and $\mathfrak{r}$ is the radical of
$\sp(1,n+1)_{\H p}$. We may write the $\mathbb{Z}$-grading
$$\sp(1,n+1)_{\H p}=\g_0+\g_1+\g_2, \quad
\g_0=\Real\oplus\spp(1)\oplus\sp(n),\quad \g_1=\H^n,\quad
\g_2=\Im\H$$ with the grading element
$1\in\Real\subset\sp(1,n+1)_{\H p}$, i.e.
$\ad_1|_{\g_\alpha}=\alpha\id_{\g_\alpha}$, $\alpha=0,1,2$.

Let $m,m_1,m_2$ be integers such that either $m+m_1+m_2= n$ or
$m+m_1+m_2\leq  n-2$ . Set the following denotation:
\begin{align*} \H^m&=\spa_\H\{e_1,...,e_m\},\\
\Im\H^{m_1}&=i\Real^{m_1}\oplus j\Real^{m_1}\oplus
k\Real^{m_1},\quad\text{\rm where}\quad\Real^{m_1}=\spa_\Real\{e_{m+1},...,e_{m+m_1}\},\\
\Co^{m_2}&=\spa_{\Real\oplus
i\Real}\{e_{m+m_1+1},...,e_{m+m_1+m_2}\}.\end{align*} Let $L'$ be
a real vector subspace of $\spa_\H\{e_{m+m_1+m_2+1},...,e_n\}$
coinciding with a $g$-orthogonal direct sum of the real spaces of
the form $$\spa_\Real\{f_1,...,f_{l}, if_1+jf_2,...,i f_{l-1}+j
f_l\},\quad l\geq 2,$$ where we fix a fragmentation of the
interval $[m+m_1+m_2+1,...,n]$ of natural numbers into a disjunct
union of subintervals of length at least 2 and $f_1,...,f_l$ are
vectors from the set $\{e_{m+m_1+m_2+1},...,e_n\}$ corresponding
to one of these subintervals.

Consider the following real vector subspace of $\H^n$:
\begin{equation}\label{Lmm1}
L=L(m,m_1,m_2,L')=\H^m\oplus\Im\H^{m_1}\oplus\Co^{m_2}\oplus
L'.\end{equation}

Assume that the decomposition \eqref{Lmm1} is $g$-orthogonal. Let
$g$ be defined by this and the following conditions:
\begin{itemize}
\item[1)] $g_{ab}=\delta_{ab}$, if $1\leq a,b\leq m$;
\item[2)] $g_{ab}=\delta_{ab}+iw_{1ab}+jw_{2ab}+kw_{3ab}$, if $m+1\leq a,b\leq m+m_1$, where
$w_1,w_2,w_3$ are skew-symmetric bilinear  forms on $\Real^{m_1}$;
\item[3)] $g_{ab}=\delta_{ab}+w_{ab}j$, if $m+m_1+1\leq a,b\leq m+m_1+m_2$, where
$w$ is a  skew-symmetric $\Co$-bilinear  form on $\Co^{m_2}$;
\item[4)] $g_{ab}=\eta_{ab}+i\Omega_{1ab}+j\Omega_{2ab}+k\Omega_{3ab}$, if $m+m_1+m_2+1\leq a,b\leq n$,
where $\eta_{ab}$ is a positive definite symmetric bilinear form
on $\spa_\Real\{e_{m+m_1+m_2+1},...,e_n\}$ and
$\Omega_1,\Omega_2,\Omega_3$ are skew-symmetric bilinear forms on
$\spa_\Real\{e_{m+m_1+m_2+1},...,e_n\}$.
\end{itemize} The above forms may be degenerate or zero.

Recall that any subalgebra $\h\subset\sp(n)$ can be decomposed as
$\h=\h'\oplus\z(\h)$, where $\h'=[\h,\h]$ is the commutant  of
$\h$ and $\z(\h)$ is the center of $\h$.

Now we state the classification result. In the first theorem we
provide a general description of possible holonomy algebras. In
the second theorem we give the precise list of all possible
holonomy algebras. It will be enough to prove only the second
theorem.

\begin{theorem}\label{Clas0} Let $(M,h)$ be a locally indecomposable pseudo-hyper-K\"ahlerian manifold
of signature $(4,4n+4)$, $n\geq 1$. If the holonomy algebra $\g$
of $(M,h)$ is not irreducible, then $\g$ is conjugated by an
element of $\SO(4,4n+4)$ to one of the following  subalgebras of
$\sp(1,n+1)_{\H p}$:
\begin{description}\item[I.]
$\g_{I}=\bar\h\zr(\mathbb{H}^{m}\oplus\mathbb{C}^{n-m}\zr\Im\H)$,
where $0\leq m\leq n$,

if $m=n$, then $\bar\h\subset\Real\oplus\sp(1)\oplus\sp(n)$ is a
subalgebra such that $\dim\pr_{\H}\bar\h\neq 1$;

if $m<n$, then $\bar\h\subset\Real\oplus i\Real\oplus\sp(m)$ is a
subalgebra such that $\dim\pr_{\Real\oplus i\Real}\bar\h\neq 1$;

\item[II.]
$\g_{II}=\h\zr(L(m,m_1,m_2,L')\zr\Im\H)$, where $\h\subset\sp(m)$
is a subalgebra.

\item[III.]
$\g_{III}=\bar\h\zr(\mathbb{H}^{m}\oplus\Im\mathbb{H}^{n-m}\zr\Im\H)$,
where $0\leq m< n$,

$\bar\h\subset\{a+\Op(aE_{n-m})|a\in\sp(1)\}\oplus\sp(m)$

and
$\pr_{\{a+\Op(aE_{n-m})|a\in\sp(1)\}}\h=\{a+\Op(aE_{n-m})|a\in\sp(1)\}$,

here $\Op(aE_{n-m})\in\sp(n-m)$ is the element with the matrix
$aE_{n-m}$.

\item[IV.] $\g_{IV}=\{A+\psi(A)|A\in\h\}\zr(\H^k\oplus V\zr\Im\H),$

here $L(m,m_1,m_2,L')=\mathbb{H}^{k}\oplus V\oplus U$ is an
$\eta$-orthogonal decomposition ($\eta=\Re g$),
$\h\subset\sp(k)$ is a subalgebra, $\psi:\h\rightarrow U$ is a
surjective linear map and $\psi\mid_{\h'}=0$.

\item[V.] $\g_{V}=\{A+\psi(A)|A\in\bar\h\}\zr(\H^{m_0}\oplus \Im\H^{n-m_0}
\zr\Im\H),$

here $0<m_0<m\leq n$,
$\bar\h\subset\{a+\Op(aE_{n-m_0})|a\in\sp(1)\}\oplus\sp(m_0)$ is a
subalgebra with\\
$\pr_{\{a+\Op(aE_{n-m_0})|a\in\sp(1)\}}\bar\h=\{a+\Op(aE_{n-m_0})|a\in\sp(1)\}$,
 $\psi:\bar\h\rightarrow
\spa_\Real\{e_{m_0+1},...,e_m\}$ is a surjective linear map with
$\psi\mid_{\bar\h'}=0$.
\end{description}
\end{theorem}

\begin{theorem}\label{Clas} Let $(M,h)$ be a locally indecomposable pseudo-hyper-K\"ahlerian manifold
of signature $(4,4n+4)$, $n\geq 1$. If the holonomy algebra $\g$
of $(M,h)$ is not irreducible, then $\g$ is conjugated by an
element of $\SO(4,4n+4)$ to one of the following  subalgebras of
$\sp(1,n+1)_{\H p}$:
\begin{description}\item[I.1.]
$\g_{I1}=\mathbb{R}\oplus\h_{0}\oplus\h\zr(\mathbb{H}^{m}\oplus\mathbb{C}^{n-m}\zr\Im\H)$,
where $0\leq m\leq n$, $\h_0\subset\spp(1)$, $\h\subset\sp(m)$ are
subalgebras, $\h_{0}=\mathbb{R}i$ or $\h_{0}=\spp(1)$. If $m<n$,
then $\h_{0}=\mathbb{R}i$.
\item[I.2.]
$\g_{I2}=\mathbb{R}\oplus\{\phi(A)+A\,|\,
A\in\h\}\zr(\mathbb{H}^{m}\oplus\mathbb{C}^{n-m}\zr\Im\H)$,\\
where $1\leq m\leq n$, $\h\subset\sp(m)$ is a subalgebra,
$\phi:\h\rightarrow\spp(1)$ is a non-zero homomorphism.

If $m<n$, then $\Im\phi=\mathbb{R}i$, $\phi\mid_{\h'}=0$. If
$m=n$, then either $\Im\phi=\mathbb{R}i$ and $\phi\mid_{\h'}=0$,
or $\Im\phi=\spp(1)$.
\item[I.3.]
$\g_{I3}=\h_{0}\oplus\{\varphi(A)+A \,
|\,A\in\h\}\zr(\mathbb{H}^{m}\oplus\mathbb{C}^{n-m}\zr\Im\H)$,\\
where $0\leq m\leq n$, $\h_0\subset\spp(1)$, $\h\subset\sp(m)$ are
subalgebras, $\varphi:\h\rightarrow\mathbb{R}$ is a linear map,
$\varphi\mid_{\h'}=0$.

If $m<n$, then $\h_{0}=\mathbb{R}i$ and $\varphi\neq 0$. If $m=n$,
then either $\h_{0}=\mathbb{R}i$ and $\varphi\neq 0$, or
$\h_{0}=\spp(1)$.
\item[I.4.]
$\g_{I4}=\{\varphi(A)+\phi(A)+A
\,|\,A\in\h\}\zr(\mathbb{H}^{m}\oplus\mathbb{C}^{n-m}\zr\Im\H)$,\\
where $0\leq m\leq n$, $\h\subset\sp(m)$ is a subalgebra,
$\varphi:\h\rightarrow\mathbb{R}$, $\phi:\h\rightarrow\spp(1)$ are
homomorphisms.

If $m<n$, then either $\varphi=\phi=0$ or $\varphi\neq 0$,
$\Im\phi=\mathbb{R}i$ and the maps $i\varphi,
\phi:\h\rightarrow\mathbb{R}i$ are not proportional,
$\varphi\mid_{\h'}=\phi\mid_{\h'}=0$.
\item[I.5.] $\g_{I5}=\Real(\alpha+i)\oplus\{\varphi(A)+A
\,|\,A\in\h\}\zr(\mathbb{H}^{m}\oplus\mathbb{C}^{n-m}\zr\Im\H)$,
where $0\leq m\leq n$, $\alpha\in\Real$, $\alpha\neq 0$,
$\h\subset\sp(m)$ is a subalgebra,
$\varphi:\h\rightarrow\mathbb{R}$ is a non-zero linear map with
$\varphi|_{\h'}=0$.
\item[II.]
$\g_{II}=\h\zr(L(m,m_1,m_2,L')\zr\Im\H)$, where $\h\subset\sp(m)$
is a subalgebra.
\item[III.1.]
$\g_{III1}=\{a+\Op(aE_{n-m})|a\in\sp(1)\}\oplus\h\zr(\H^m\oplus\im
\H^{n-m}\zr\Im\H)$, where $n-m\geq 1$, $\h\subset\sp(m)$ is a
subalgebra, and $\Op(aE_{n-m})\in\sp(n-m)$ is the element with the
matrix $aE_{n-m}$.
\item[III.2.]
$\g_{III2}=\{\phi(A)+A+\Op(\phi(A)E_{n-m})|A\in\h\}\zr(\H^m\oplus\im
\H^{n-m}\zr\Im\H)$,\\ where $n-m\geq 1$, $\h\subset\sp(m)$ is a
subalgebra, and $\phi:\h\to\sp(1)$ is a surjective homomorphism.
\item[IV.] $\g_{IV}=\{A+\psi(A)|A\in\h\}\zr(\H^k\oplus
V\zr\Im\H),$\\
here $L(m,m_1,m_2,L')=\mathbb{H}^{k}\oplus V\oplus U$ is an
$\eta$-orthogonal decomposition ($\eta=\Re g$),
$\h\subset\sp(k)$ is a subalgebra, $\psi:\h\rightarrow U$ is a
surjective linear map and $\psi\mid_{\h'}=0$.

\item[V.1.] $\g_{V1}=\{a+\Op(aE_{n-m_0})|a\in\sp(1)\}\oplus\{A+\psi(A)|A\in\h\}\zr(\H^{m_0}\oplus \Im\H^{n-m_0}
\zr\Im\H),$\\ here $0<m_0<m\leq n$, $\h\subset\sp(m_0)$ is a
subalgebra, $\psi:\h\rightarrow \spa_\Real\{e_{m_0+1},...,e_m\}$
is a surjective linear map with $\psi\mid_{\h'}=0$.

\item[V.2.] $\g_{V2}=\{a+\chi(a)+\Op(aE_{n-m_0})|a\in\sp(1)\}\oplus\{A+\psi(A)|A\in\h\}
\zr(\H^{m_0}\oplus \Im\H^{n-m_0} \zr\Im\H),$\\ where the dates are
the same as for $\g_{V1}$ and in addition $\chi:\sp(1)\to\sp(m_0)$
is an injective homomorphism such that $\chi(\sp(1))$ commutes
with $\h$.
\end{description}
Conversely, all these algebras are Berger algebras.
\end{theorem}

{\bf Another description of $\sp(1,n+1)_{\H p}$}

For convenience we give another description of the Lie algebra
$\sp(1,n+1)_{\H p}$. One usually identifies the Lie algebra
$\so(r,s)$ with the space of bivectors $\Lambda^2\Real^{r,s}$ in
such a way that $(X\wedge Y)Z= \eta(Z,X)Y- \eta(Z,Y)X.$ Having the
pseudo-quaternionic-Hermitian metric $g$ on $\H^{1,n+1}$ we may
put $$(X\wedge_g Y)Z= g(Z,X)Y- g(Z,Y)X,$$ then
$$X\wedge_g Y=\sum_{s=0}^3 I_rX\wedge  I_sY,\quad (a X)\wedge_g Y=X\wedge_g (\bar aY),$$
and $X\wedge_g Y\in\sp(1,n+1)$. We get the identification
$$\sp(1,n+1)\simeq\Lambda_g^2\H^{1,n+1}=\{X\wedge_gY|\, X,Y\in\H^{1,n+1}\}.$$
The element $(a,A,X,b)\in\sp(1,n+1)_{\H p}$ corresponds to
$$-(ap)\wedge_g q+A+p\wedge_g X+p\wedge_g \left(\frac{b}{2}p\right),\quad A\in\sp(n)\simeq \Lambda_g^2\H^{n}.$$
Let us rewrite the Proposition \ref{propR} given below in these
notations.

\begin{prop}\label{propRA} Any $R\in\mathcal{R}(\sp(1,n+1)_{\H p})$ is uniquely
defined by elements $C_{1},C_{2}\in\H$,
$A_{1},A_{2},A_{3}\in\sp(n)$, $S_{1},S_{2}\in\H^n$,
$R'\in\R(\sp(n))$, $P\in\P(\sp(n))$, $d_1,...,d_5\in\Real$ in the
following way:
\begin{align*} R(I_sp,q)&=p\wedge_g\left(\frac{1}{2} B_{s}p\right) ,\qquad
R(q,X)=P(X)+p\wedge_g T(X)+p\wedge_g\left(\frac{1}{2}\theta(X)p\right),\\
R(X,Y)&=R'(X,Y)+p\wedge_g
(P(Y)X-P(X)Y)+p\wedge_g\left(\frac{1}{2}\Big(g(Y,T(X))-g(X,T(Y))\Big)p\right),\\
R(q,I_sq)&=-(C_{s}p)\wedge_g  q+ A_{s}+p\wedge_g S_{s}+p\wedge_g \left(\frac{1}{2}D_{s}q\right),\\
R(p,I_sp)&=R(p,X)=0,\qquad X,Y\in\H^{n},\end{align*} where
\begin{align*} C_{3}&=C_{2}i-C_{1}j,\quad T=-\frac{1}{2}(I_1A_{1}+I_2A_{2}+I_3A_{3}),\quad S_{3}=jS_{1}-iS_{2},\\
 D_{1}&=d_{1}i+d_{2}j+d_{3}k,\quad
D_{2}=d_{2}i+d_{4}j+d_{5}k,\quad D_{3}=jD_{1}-iD_{2}, \\
B_{s}&=\frac{1}{2}(I_1I_sC_{1}+I_2I_sC_{2}+I_3I_sC_{3}),\quad
\theta(X)=\frac{1}{2}(I_1g(X,S_{1})+I_2g(X,S_{2})+I_3g(X,S_{3})).
\end{align*}
\end{prop}
The other values of $R$ can be found using equality \eqref{RI}
given below. Although this version of Proposition \ref{propR} is
not so complicated, the form of Proposition \ref{propR} is more
convenient for the proof of Theorem \ref{Clas}.

\section{Proof of Theorem \ref{Clas}}\label{D}

Since $\g$ is weakly irreducible and not irreducible, $\g$
preserves a degenerate vector subspace $V\subset\Real^{4,4n+4}$.
Let $V_1=V\cap V^\bot$, then $V_1$ is isotropic and $\dim V_1\leq
4$. Let $V_2=V_1^\bot\cap I_1 V_1^\bot$. Clearly, $V_2\neq 0 $ and
it is degenerate, $\g$-invariant and $I_1$-invariant. Then
$V_3=V_2\cap V_2^\bot$ is isotropic, $\g$-invariant and
$I_1$-invariant. Starting with $V_3$ in the same way it can be
shown that $\g$ preserves an isotropic $I_1,I_2$-invariant
subspace $W\subset\Real^{4,4n+4}$, then $W$ is also
$I_3$-invariant and it has dimension 4. Consequently,
$\g\subset\sp(1,n+1)_{\H p}$.

The proof of the Theorem will consist of several parts.

\subsection{The structure of the space $\mathcal{R}(\sp(1,n+1)_{\H p})$ }

Let us find the space of curvature tensors
$\mathcal{R}(\sp(1,n+1)_{\H p})$ for the Lie algebra
$\sp(1,n+1)_{\H p}$.

Using the form $\eta$, the Lie algebra $\so(4,4n+4)$ can be
identified with the space
$$\wedge^{2}\mathbb{R}^{4,4n+4}=\spa\{u\wedge v=u\otimes v-v\otimes u
|u,v\in\mathbb{R}^{4,4n+4}\}$$ in such a way that $(u\wedge
v)w=\eta(u,w)v-\eta(v,w)u$ for all $u,v,w\in \mathbb{R}^{4,4n+4}$.
 One can check that the element
$\Op\left(\begin{array}{ccc}a&-(G\bar{X})^{t}&b\\0&\Mat_A&X\\0&0&-\bar
a\end{array}\right)\in\sp(1,n+1)_{\H p}$ corresponds to the
bivector
\begin{align*}
&\Big(-a_0(p\wedge q+ip\wedge iq+ jp\wedge jq+ kp\wedge
kq)+a_1(p\wedge iq-ip\wedge q+ kp\wedge jq- jp\wedge kq)\\ &+
a_2(-jp\wedge q+p\wedge jq- kp\wedge iq+ ip\wedge
kq)+a_3(-kp\wedge q+jp\wedge iq-ip\wedge jq+ p\wedge kq)\Big)\\
&+ A+\Big((p\wedge X_{0}+ip\wedge iX_{0}+jp\wedge jX_{0}+kp\wedge
kX_{0})+(p\wedge iX_{1}-ip\wedge X_{1}-jp\wedge kX_{1}+kp\wedge
jX_{1})\\ &+(p\wedge jX_{2}+ip\wedge kX_{2}-jp\wedge
X_{2}-kp\wedge iX_{2})+(p\wedge kX_{3}-ip\wedge jX_{3}+jp\wedge
iX_{3}-kp\wedge X_{3})\Big)\\&+ \Big(b_1(p\wedge ip-jp\wedge
kp)+b_2(p\wedge jp+ip\wedge kp)+b_3(p\wedge kp-ip\wedge
jp)\Big),\end{align*} where $X=X_0+iX_1+jX_2+kX_3$,
$X_0,...,X_3\in\Real^n=\spa_\Real\{e_1,...,e_n\}\subset\Real^{4n}\simeq\H^n,$
$A\in\sp(n)\subset\so(4n)\simeq\wedge^2\Real^{4n}$.

Let $\g\subset\sp(1,n+1)_{\H p}$, $R\in\mathcal{R}(\g)$. It is
known that any $R\in\mathcal{R}(\g)$ satisfies
\begin{equation}\label{(*)}
\eta(R(u, v)z,w)=\eta(R(z, w)u,v)
\end{equation} for all $u,v,w,z\in\Real^{4,4n+4}$. Using this property, we get the following:
\begin{multline*}\eta(R(I_{\alpha}u,v)z,w)=\eta(R(z,w)I_{\alpha}u,v)=
\eta(I_{\alpha}R(z,w)u,v)=-\eta(R(z,w)u,I_{\alpha}v)=-\eta(R(u,I_{\alpha}v)z,w),\end{multline*}
i.e. for any  $1\leq \alpha\leq 3$ and $u,v\in\Real^{4,4n+4}$,
\begin{equation}\label{RI}R(I_\alpha u,v)= -R(u,I_\alpha v)\end{equation} holds.
Hence,
\begin{equation}R(xu,v)=R(u,\bar xv)\end{equation}
for all $x\in \H$ and $u,v\in\Real^{4,4n+4}$.

The metric $\eta$ defines the metric $\eta\wedge\eta$ on
$\wedge^{2}\mathbb{R}^{4,4n+4}$. Using the above identification,
$R$ can be considered as the map
$R:\wedge^{2}\mathbb{R}^{4,4n+4}\to\g\subset\so(4,4n+4)\simeq\wedge^2\mathbb{R}^{4,4n+4}$.
From \eqref{(*)}, we obtain
\begin{equation}\label{symR}\eta\wedge\eta(R(u\wedge v), z\wedge
w)=\eta\wedge\eta(R(z\wedge w), u\wedge v)\end{equation} for all
$u,v,z,w\in\mathbb{R}^{4,4n+4}$. This shows that $R$ is a
symmetric linear map. Consequently $R$ is zero on the orthogonal
complement to $\g$ in $\wedge^{2}\mathbb{R}^{4,4n+4}$. In
particular, the vectors $q\wedge iq+jq\wedge kq,$ $q\wedge
jq-iq\wedge kq,$ $q\wedge kq+iq\wedge jq$, $I_{r}p\wedge I_{s}p$,
$I_{r}p\wedge X$, where $X\in\mathbb{R}^{4n}$, are contained in
the orthogonal complement to $\sp(1,n+1)_{\H p}$.
Hence,\begin{equation}\label{(***)} R(q, iq)=-R(jq, kq),\, R(q,
jq)=R(iq, kq),\, R(q, kq)=-R(iq, jq),\ R(I_{r}p, I_{s}p)=R(I_rp,
X)=0.
\end{equation}

For a subalgebra $\h\subset\so(m)$ define the space
$$\P(\h)=\{P\in (\Real^m)^*\otimes
\h|\eta(P(x)y,z)+\eta(P(y)z,x)+\eta(P(z)x,y)=0\text{ for all }
x,y,z\in \Real^m\},$$ where $\eta$ is the scalar product on
$\Real^m$. This space is studied in \cite{Leistner,GalP}.

\begin{prop}\label{propR} Any $R\in\mathcal{R}(\sp(1,n+1)_{\H p})$, $n\geq 1$ is uniquely
defined by elements $C_{01},C_{02}\in\H$,
$A_{01},A_{02},A_{03}\in\sp(n)$, $S_{01},S_{02}\in\H^n$,
$R'\in\R(\sp(n))$, $P_0\in\P(\sp(n))$, $d_1,...,d_5\in\Real$ in
the following way:
\begin{align*} R(I_rp,I_sq)&=(0,0,0,B_{rs}),\qquad\qquad\qquad\qquad\qquad
R(I_sq,X)=(0,P_s(X),T_s(X),\theta_s(X)),\\
R(X,Y)&=(0, R'(X,Y),Q(X,Y),\tau(X,Y)),\qquad
R(I_rq,I_sq)=(C_{rs}, A_{rs},S_{rs},D_{rs}),\\
R(I_rp,I_sp)&=R(I_rp,X)=0,\qquad X,Y\in\Real^{4n},\end{align*}
where \begin{align}\label{p1C} C_{03}&=C_{02}i-C_{01}j,\,
C_{12}=-C_{03},\,C_{13}=C_{02},\,C_{23}=-C_{01},\,C_{rs}=C_{0r}I_{s}-C_{0s}I_{r},\,r,s\neq 0,\\
\label{p1D1} D_{01}&=d_{1}i+d_{2}j+d_{3}k,\quad
D_{02}=d_{2}i+d_{4}j+d_{5}k,\quad D_{03}=jD_{01}-iD_{02}, \\
\label{p1D2} D_{23}&=-D_{01},\quad D_{13}=D_{02},\quad
D_{12}=-D_{03},\quad
D_{rs}=I_rD_{0s}-I_sD_{0r}, \\
\label{p1A} A_{23}&=-A_{01},\quad A_{13}=A_{02},\quad
A_{12}=-A_{03},\\\label{p1To}
T_0&=-\frac{1}{2}(I_1A_{01}+I_2A_{02}+I_3A_{03}),\quad
T_s=I_sT_0-A_{0s}=-T_0I_s,\quad s\neq 0,\\
\label{p1P} P_s&=-P_0\circ I_s,\quad s\neq 0,\quad
Q(X,Y)=P_0(Y)X-P_0(X)Y,\\
\label{p1B} B_{rs}&=I_rC_{0s}+I_sB_{r0},\quad
B_{r0}=\frac{1}{2}(I_1I_rC_{01}+I_2I_rC_{02}+I_3I_rC_{03}),\\
\label{p1tau}
\tau(X,Y)&=g(Y,T_0(X))-g(X,T_0(Y)),\\
\label{p1S} S_{03}&=jS_{01}-iS_{02},\, S_{23}=-S_{01},\,
S_{13}=S_{02},\, S_{12}=-S_{03},\,
S_{rs}=I_rS_{0s}-I_sS_{0r}, \\
\label{p1theta0}\theta_0(X)&=\frac{1}{2}(I_1g(X,S_{01})+I_2g(X,S_{02})+I_3g(X,S_{03})),\\
\label{p1thetas}
\theta_s(X)&=g(X,S_{0s})+I_s\theta_0(X)=-\theta_0(I_sX), \quad
s\neq 0,
\end{align} where  $X,Y\in\H^n$. Moreover, it holds
\begin{multline}\label{LP}\eta(Q(Y,Z),X)=\eta(P_0(X)Y,Z),\quad\eta(I_r\tau(X,Y)p,I_sq)=\eta(A_{rs}X,Y),\\
\eta(I_r\theta_s(x)p,I_tq)=\eta(I_sS_{rt},X),\quad \eta
(I_{t}B_{rs}p,I_{t_{1}}q)=\eta
(I_{r}C_{tt_{1}}p,I_{s}q),\end{multline} where
$X,Y,Z\in\Real^{4n}$.
\end{prop}

{\bf Proof.} Let $R\in\mathcal{R}(\sp(1,n+1)_{\H p})$.

The equality \eqref{(***)} shows that
$R(I_{r}p,I_{s}p)=R(I_rp,X)=0$. We may write
\begin{align*} R(I_rp,I_sq)&=(\lambda_{rs}, F_{rs},X_{rs},B_{rs}),\qquad\qquad\qquad
R(I_sq,X)=(\mu_s(X),P_s(X),T_s(X),\theta_s(X)),\\
R(X,Y)&=(\sigma(X,Y), R'(X,Y),Q(X,Y),\tau(X,Y)),\qquad
R(I_rq,I_sq)=(C_{rs}, A_{rs},S_{rs},D_{rs}).\end{align*} Now we
find the conditions that satisfy the obtained elements. By the
Bianchi identity, $$R(I_rp,I_sq)X+R(I_sq,X)I_rp+R(X,I_rp)I_sq=0.$$
Using the equality $R(X,I_rp)=0$ and taking the projection on
$\H^{n}$, we get $ F_{rs}X=0$, i.e. $F_{rs}=0$. Using \eqref{(*)},
we obtain
$$\eta(R(X,Y)I_rp,I_sq)=\eta(R(I_rp,I_sq)X,Y)=0,$$ hence
$\sigma(X,Y)=0$. As in \cite[Proposition 1]{44}, one can prove
that $\lambda_{rs}=0$ and the equalities for
$C_{rs},D_{rs},B_{rs}$. Writing down the Bianchi identity for the
vectors $I_rp,I_sq,I_tq$ and taking the projection on $\H^n $, we
have
$$I_tX_{rs}=I_sX_{rt}.$$
Hence, $X_{rs}=I_sX_{r0}$. Substituting this back to the above
equation, we get $I_tI_sX_{r0}=I_sI_tX_{r0}$. Taking $t=1$, $s=2$,
we obtain $X_{r0}=0$. This shows that $X_{rs}=0$. From this and
\eqref{(*)} it follows that $\mu_s=0$.

Writing down the Bianchi identity for the vectors
$I_{r}q,I_{s}q,I_{t}q$ and taking the projection on $\H^n $, we
get $$I_{t}S_{rs}+I_{r}S_{st}+I_{s}S_{tr}=0.$$ Taking $t=0$, we
have
$$S_{rs}=I_rS_{0s}-I_sS_{0r}.$$ Substituting this to the initial
equality and taking $t=1, r=2, s=3$, we obtain
$$S_{03}=jS_{01}-iS_{02}.$$ Note that
$R(X,Y)Z=R'(X,Y)Z-g(Z,Q(X,Y))p$. The Bianchi identity written for
the vectors $X,Y,Z$ implies $R'\in\R(\sp(n))$. Moreover,
\begin{equation}g(Z,Q(X,Y))+g(X,Q(Y,Z))+g(Y,Q(Z,X))=0.\end{equation}
Hence,
\begin{equation}\label{(L)}\eta(Q(X,Y),Z)+\eta(Q(Y,Z),X)+\eta(Q(Z,X),Y)=0.\end{equation}
From the Bianchi identity written for the vectors $q,X,Y$ it
follows that $$P_0(X)Y+Q(X,Y)-P_0(Y)X=0.$$ This and \eqref{(L)}
imply $P_0\in\P(\h)$. Using \eqref{RI}, we get
$R(I_sq,X)=-R(q,I_sX)$, hence $$P_s(X)=-P_0(I_sX),\quad
T_s(X)=-T_0(I_sX),\quad \theta_s(X)=-\theta_0(I_sX),\quad  s\neq
0.$$ The Bianchi identity written for the vectors $I_rq,I_sq,X$
implies
\begin{equation}\label{AT}
A_{rs}X+I_rT_s(X)-I_sT_r(X)=0.\end{equation} Taking $r=0$, we get
$T_s(X)=I_sT_0(X)-A_{0s}(X)$. Substituting this to \eqref{AT} and
taking $r=1$, $s=2$, we obtain
$$T_0=-\frac{1}{2}(I_1A_{01}+I_2A_{02}+I_3A_{03}).$$ Writing down
the Bianchi identity for the vectors $I_rq,I_sq,X$ and taking the
projection on $\H p$, we get
$$-g(X,S_{rs})+I_r\theta_s(X)-I_r\theta_r(X)=0.$$
Using this it is easy to get \eqref{p1theta0} and
\eqref{p1thetas}. The Bianchi identity applied to $X,Y,q$ implies
the equality for $\tau(X,Y)$ from \eqref{p1tau}. We have proved
that any $R\in\R(\sp(1,n+1)_{\H p})$ satisfies the conditions of
the proposition.

Conversely, it can be checked that any element $R$ satisfying  the
conditions of the proposition belongs to $\R(\sp(1,n+1)_{\H p})$.
\hfill $\square$

Denote by $\sp(1,1)_{\H p}$ the subalgebra of $\sp(1,n+1)_{\H p}$
that annihilates $\H^n\subset\H^{1,n+1}$. The space
$\R(\sp(1,1)_{\H p})$ is found in \cite[Proposition 1]{44}. Note
that any $R$ given by elements $C_{rs}$, $B_{rs}$, $D_{rs}$ and
such that all the  rest elements are zero belongs to
$\R(\sp(1,1)_{\H p})$. In particular, we get

\begin{lemma}\label{L1}
Any subalgebra   $\g\subset\sp(1,n+1)_{\H p}$ such that
$\dim_\Real\pr_{\H}\g=1$  is not a Berger algebra. \qed
\end{lemma}

\subsection{The algebras listed in the statement of the
theorem \ref{Clas} are Berger algebras}

Here we prove that the algebras listed in the statement of the
theorem \ref{Clas} are Berger algebras.

Let $\g=\g_{I1}$. If $m=n$ and $\h_0=\spp(1)$, then any
$R\in\R(\g)$ is given as in Proposition \ref{propR} with
$A_{01},A_{02},A_{03}\in\h$. Since the elements $C_{01}\in\H$,
$A_{01}\in\h$, $S_{01}\in\H^n$, $D_{01}\in\Im\H$ can be chosen in
arbitrary way, $\g$ is a Berger algebra. If $m=n$ and $\h_0=\Real
i$, then from \cite[Section 4]{44} it follows that in addition to
the above case $C_{01}=0$ and $C_{02}\in\Real\oplus\Real i$ is
arbitrary, hence $\g$ is a Berger algebra. Suppose that $m<n$.
Then $\h_0=\Real i$. Each $S_{rs}$ can be written as
$S_{rs}=S'_{rs}+S''_{rs}$, where $S'_{rs}\in\H^m$ and
$S''_{rs}\in\Co^{n-m}$. Then
$S''_{01},S''_{02},S''_{03}\in\Co^{n-m}$. The condition
$S''_{03}=jS''_{01}-iS''_{02}$ implies $S''_{01}=0$, on the other
hand, $S''_{02}\in\Co^{n-m}$ is arbitrary. This shows that $\g$ is
a Berger algebra.

Let $\g=\g_{I2}$. Suppose that $\im \phi=\Real i$. Let
$R\in\R(\g)$. From the above example we get that $C_{01}=0$. In
addition, $C_{02}=c_1+\phi(A_{02})$, where $c_1\in\Real$. Hence,
$C_{03}=\phi(A_{02})i+c_1i$. This shows that $c_1=-i\phi(A_{03})$.
Consequently $\g$ is a Berger algebra. The case $\im\phi=\sp(1)$
will follow from this and the next case.

Let $\g=\g_{I4}$ and $R\in\R(\g)$. Let
$\phi=i\phi_1+j\phi_2+k\phi_3$, where the maps
$\phi_1,\phi_2,\phi_3$ take values in $\Real$. Then
$$C_{rs}=\varphi(A_{rs})+\phi_1(A_{rs})i+
\phi_2(A_{rs})j+\phi_3(A_{rs})k.$$ The condition
$C_{03}=C_{02}i-C_{01}j$  is equivalent to the equalities
\begin{align*}\phi_{2}(A_{01})-\phi_{1}(A_{02})&=\varphi(A_{03}),\quad
\varphi(A_{02})+\phi_{3}(A_{01})=\phi_{1}(A_{03}),\\
\phi_{3}(A_{02})-\varphi(A_{01})&=\phi_{2}(A_{03}),\quad
-\phi_{2}(A_{02})-\phi_{1}(A_{01})=\phi_{3}(A_{03}).\end{align*}
It is not hard to see that these conditions can be satisfied
taking appropriate $A_{01},A_{02},A_{03}$. For example, if
$\im\phi=\Real i$, then there exists a decomposition
$\h=\h_1\oplus\h_2\oplus\h_3$ such that
$\ker\varphi={\h_2\oplus\h_3}$ and $\ker\phi={\h_1\oplus\h_3}$. In
this case, it is enough to take $A_{01}\in\h_3$, $A_{02}\in\h_2$,
$A_{03}\in\h_1$ such that $\varphi(A_{03})=-\phi_1(A_{02})=1$.
This shows that $\g$ is a Berger algebra.

The other Lie algebras from the statement of the theorem can be
considered in the same way. For $\g_{III2}$ and $\g_{IV}$ note the
following. Obviously, $L'$ satisfies condition \eqref{cond3} given
below. Let $X,Y,jX-iY\in L'$ and $S_{01}=X,$ $S_{02}=Y$, then
$S_{03}=jS_{01}-iS_{02}=jX-iY\in L'$, hence $L'$ is spanned by
$S_{rs}$.

\subsection{Weakly-irreducible subalgebras of $\sp(1,n+1)$ and
real vector subspaces in $\H^n$}\label{w}

Now we review  the classification of  weakly irreducible
subalgebras $\g\subset\sp(1,n+1)_{\H p}$ obtained in \cite{B1} and
make several corrections.

In \cite{B1} was constructed a homomorphism $f:\sp(1,n+1)_{\H
p}\to\simil\H^n$, where $ \simil
\H^n=\Real\oplus(\sp(1)\oplus\sp(n))\zr\H^n$ is the Lie algebra of
the group $\Sim\H^n$ of similarity transformations of $\H^n$. The
homomorphism $f$ is surjective with the kernel $\Im\H$ and it is
given as $f(a_0+a_1,A,X,b)=(a_0,-a_1+A,X)$, where $a_0\in\Real$,
$a_1\in\sp(1)$. Let $\g\subset\sp(1,n+1)_{\H p}$ be a weakly
irreducible subalgebra and $L=\pr_{\H^n}\g\subset \H^n$. It is
shown that  $\spa_{\H} L=\H^n$ and the connected subgroup of
$\Sim\H^n $ corresponding to $f(\g)\subset \simil \H^n $ preserves
$L$ and acts on it transitively.  It was stated that there exists
a $g$-orthogonal decomposition $L=L_1\oplus L_2\oplus L_3$ such
that $L_1\subset\H^n$ is a quaternionic subspace, $L_2\subset
\H^n$ is a real subspace such that $iL_2=L_2$, $L_2\cap jL_2=0$,
$\spa_\H L_2=L_2\oplus jL_2$, and $L_3\subset \H^n$ is a real
subspace such that $L\cap I_sL=0$, $1\leq s\leq 3$. Let $n=1$ and
$e\in\H$ be a non-zero vector, then the vector subspace $L=i\Real
e\oplus j\Real e\oplus k\Real e$ is missing. Next, it was stated
that $\dim_\Real L_3=\dim_\H\spa_\H L_3$. Let $n=2$, $e_1,e_2$ a
basis of $\H^2$, then the space
$L=\spa_\Real\{e_1,e_2,je_1+ie_2\}$ is missing. Finally it was
stated that there exists a $g$-orthonormal basis of $\H^n$
consisting of vectors from $L_1$, $L_2$ and $L_3$, which is also
not true in general.

Let $L\subset\H^n$ be a real subspace such that $\spa_\H L=\H^n$.
Put $L_1=L\cap iL\cap jL\cap kL$, i.e. $L_1$ is the biggest
quaternionic vector subspace in $L$. Let $L_2=\{X\in
L|g(X,L_1)=0\}$,
 then $L=L_1\oplus L_2$ and $L_2\cap iL_2\cap jL_2\cap kL_2=0.$
Note that the subspaces $L_2\cap iL_2\cap jL_2$, $L_2\cap iL_2\cap
kL_2$, $L_2\cap jL_2\cap kL_2$, $iL_2\cap jL_2\cap kL_2$ can be
taken to each other by $i,j,k$, and the intersection of any two of
these subspaces are zero. Let  $L_3$ be the direct sum of these
subspaces. Let $L_4=\{X\in L_2|g(X,L_3)=0\}$ and $L_5=\{Y\in
L_2|g(Y,L_4)=0\}$. Then $L_2=L_5\oplus L_4$ and
$$L_5=L_3\cap L_2=iU\oplus jU\oplus kU=\sp(1)\cdot U,\quad \sp(1)=\spa_\Real\{i,j,k\}, \quad
U=iL_2\cap jL_2\cap kL_2.$$ By the construction, it holds $I_r
L_4\cap I_s L_4\cap I_t L_4=0$, if $r,s,t$ are pairwise different.
We see that $L_5$ is the biggest subspace of $L_2$ of the form
$\sp(1)\cdot V$, where $V\subset \spa_\H L_2$ is a real subspace.
In particular, this shows that the definition of $L_5$ does not
depend on the choice of the generators $I_1,I_2,I_3$ of
$\sp(1)=\spa_\Real\{I_1,I_2,I_3\}$. Let $h_1\in\sp(1)$ be an
element with $h_1^2=-\id$ (any non-zero element of $\sp(1)$ is
proportional to such element). Let $L^1_4=h_1L_4\cap L_4$. Suppose
that $L_4^1\neq 0$, then there is a $g$-orthogonal decomposition
$L_4=L_4^1\oplus L_4'$ and it holds $h_1L_4'\cap L_4'$. Taking
other $h_2\in\sp(1)$, we may decompose $L_4'$. Clearly, this
process is finite and we will get a $g$-orthogonal decomposition
$$L_4=L_4^1\oplus\cdots\oplus L_4^l\oplus L'$$
such that each $L_4^\alpha$ is $h_\alpha$-invariant for some
$h_\alpha\in\sp(1)$ with $h_\alpha^2=-\id$ and $h L_4^\alpha\cap
L_4^\alpha=0$ if $h\in\sp(1)$ is not proportional to $h_\alpha$.
Next, $h L'\cap L'=0$ for any non-zero $h\in \sp(1)$. Now we
change the quaternionic structure
$\sp(1)=\spa_\Real\{I_1,I_2,I_3\}$ on $\Real^{4n}$ to another
quaternionic structure $\widetilde{\sp(1)}=\spa_\Real\{\tilde
I_1,\tilde I_2,\tilde I_3\}$ such that $\tilde I_1|_{\spa_\H
L_4^\alpha}=h_\alpha$. This means that we consider subalgebras of
$\sp(1,n+1)$ up to a conjugancy  by elements of $\SO(4,4n+4)$.
After such change we get $$L=L_1\oplus L_5\oplus L_4^1\oplus L',$$
where $L_1,$ $L_5$, $L'$ satisfy the same properties as above and
$L_4^1$ is $I_1$-invariant.
 Let
$e_1,..,e_m$ be a $g$-orthonormal basis of $L_1\simeq\H^m$. Let
$\Real^{m_1}=iL_2\cap jL_2\cap kL_2$ and let
$\{e_{m+1},...,e_{m+m_1}\}$ be an $\eta$-orthonormal basis of
$\Real^{m_1}$. Obviously, $L_5= i\Real^{m_1}\oplus
j\Real^{m_1}\oplus k\Real^{m_1}.$ Let $X,Y\in \spa_\H L'_4$. Note
that the equality $$h(X,Y)=\eta(X,Y)+i\eta(X,I_1Y)$$ defines a
Hermitian metric on the complex space $\spa_\H L'_4$. It holds
$$g(X,Y)=h(X,Y)+h(X,I_2Y)j.$$
Let $e_{m+m_1+1},...,e_{m+m_1+m_2}$ be an $h$-orthonormal basis of
the complex space $L'_4$. Let $w(X,Y)=h(X,I_2Y)$, then the
restriction of $w$ to $L'_4$ is a $\Co$-linear skew-symmetric
bilinear  form.

To describe the structure of $L'$ we use results from \cite{DR78},
where all real subspaces $V$ of quaternionic vector spaces $U$ are
found. First  a pair $(V,U)$ of such spaces  is called
indecomposable if there are no pairs $(V_1,U_1)$, $(V_2,U_2)$ such
that $U=U_1\oplus  U_2$ and $V=V_1\oplus V_2$. In our case, $L'$
may be decomposed into a $g$-orthogonal direct sum of real spaces
$V$ such that the pair $(V,\spa_\H V)$ is indecomposable. By our
construction, it is enough to consider pairs $(V,\spa_\H V)$ such
that $hV\cap V=0$ for any $h\in\sp(1)$. Then we get only the
following two possibilities:
\begin{multline*}V=A(2l-1)=\spa_\Real\{f_1,...,f_{l-1},f_{l+1},...,f_{2l-1},\\
if_1+jf_2,...,i f_{l-1}+j f_l,\,\, f_l+if_{l+1},\,\,
jf_{l+1}+if_{l+2},...,jf_{2l-2}+if_{2l-1}\},\end{multline*} where
$f_1,...,f_{2l-1}$ ($l\geq 2$) is a basis of $\H^{2l-1}$, and
$$V=B(l)=\spa_\Real\{f_1,...,f_{l},\,\,
if_1+jf_2,...,i f_{l-1}+j f_l\},$$ where $f_1,...,f_{l}$ ($l\geq
1$) is a basis of $\H^l$.

We get that $L$ is given by
\begin{equation}\label{decLvdokve}
L=\H^m\oplus\Im\H^{m_1}\oplus\Co^{m_2}\oplus L',\end{equation}

i.e. as in Section \ref{REZ}, but at the moment $L'$ is a
$g$-orthogonal direct sum of vector spaces of the form $A(2l-1)$
and $B(l)$.

Summing the above arguments and the results from \cite{44}, we
obtain the following theorem.

\begin{theorem}\label{ThRes31}  Let $n\geq1$. Any weakly irreducible subalgebra of $\sp(1,n+1)_{\mathbb{H}{p}}$ that
contains the ideal $\Im\H$ is conjugated by an element of
$\SO(4,4n+4)$ to one of the following subalgebras:
\begin{description}

\item{{\bf Type $\Real$.}} $\g=\Real\oplus\bar{\h}\ltimes (L\ltimes\Im\H)$,

\item{{\bf Type $\varphi$.}} $\g=\{\varphi(A)+A|\,\,A\in\bar{\h}\}\ltimes (L\ltimes\Im\H)$,

\item{{\bf Type $\psi$.}} $\g=\{A+\psi(A)|\,\,A\in\bar{\h}\}\ltimes (W\ltimes\Im\H)$,
\end{description}
where $L$ is as in \eqref{decLvdokve},
$\bar{\h}\subset\sp(1)\oplus\sp(n)$ is a subalgebra such that
$\tilde\h=\{(-a,A)|(a,A)\in\bar\h\}\subset \sp(1)\oplus\sp(n)$
preserves $L$; $\varphi:\bar{\h}\to\Real$ is a linear map with
$\varphi|_{\bar{\h}'}=0$; for the last algebra $L=W\oplus U$ is an
orthogonal decomposition, the Lie algebra $\tilde\h$ annihilates
$U$, and $\psi:\bar{\h}\to U$ is a surjective linear map with
$\psi|_{\bar{\h}'}=0$.
\end{theorem}
Note that $\bar\h=\pr_{\sp(1)\oplus\sp(n)}\g$ and
$\tilde\h=\pr_{\sp(1)\oplus\sp(n)}f(\g)$.

\subsection{Classification of the Berger algebras containing $\im\H$}

Let $\g\subset\sp(1,n+1)_{\H p}$ be a weakly irreducible
subalgebra, $\bar\h=\pr_{\sp(1)\oplus\sp(n)}\g$, $\tilde
h=\{(-a,A)|(a,A)\in\bar\h\}$ and $L=\pr_{\H^n}\g$. Then
$L\subset\H^n$ is a subspace of the form
$L=\H^m\oplus\Im\H^{m_1}\oplus\Co^{m_2}\oplus L'$ (see Section
\ref{w}), and $\tilde\h$ preserves $L$. In particular, $\tilde\h$
is contained in the intersection
\begin{equation}\label{int1}
\sp(1)\oplus\sp(n)\cap\so(L)\oplus\so(L^{\bot_\eta}).\end{equation}

\begin{lemma} \label{LL'} Let $\g\subset\sp(1,n+1)_{\H p}$ be a weakly irreducible
Berger subalgebra. Then the following holds:
\begin{description}\item[1)] If $L'\neq 0$, then
$\bar\h\subset\sp(m)$, $\pr_\H\g=0$ and $L'$ is a $g$-orthogonal
sum of the spaces of type $B(l)$ with $l\geq 2$.
 \item[2)] Suppose that $L'=0$, then
 \begin{description}\item[2.a)] if $m_1\neq 0$ and $m_2= 0$, i.e. $L=\H^m\oplus\im \H^{m_1}$, $m+m_1=n$,
then $\pr_\Real\g=0$, \\ $\bar h\subset
\{a+\Op(aE_{m_1})|\,a\in\sp(1)\}\oplus\sp(m)$, the projection of
$\bar\h$ to $\{a+\Op(aE_{m_1})|\\ a\in\sp(1)\}$ is either trivial
or coincides  with $\{a+\Op(aE_{m_1})|\,a\in\sp(1)\}$;
\item[2.b)] if $m_1=0$ and $m_2\neq 0$, i.e. $L=\H^m\oplus\Co^{m_2}$, $m+m_2=n$, then
$\bar h\subset \Real i\oplus\sp(m)$;
\item[2.c)] if $m_1\neq 0$ and $m_2\neq 0$, i.e.
$L=\H^m\oplus\im \H^{m_1}\oplus\Co^{m_2}$, $m+m_1+m_2=n$, then
$\bar h\subset\sp(m)$ and $\pr_\H\g=0$.
\end{description}\end{description}
\end{lemma}

{\bf Proof.} {\bf 1)} Suppose that $L'\neq 0$. Let $R\in\R(\g)$ be
a tensor given as in Proposition \ref{propR}. Then,
$$\pr_{\H\oplus\sp(n)}R(q,I_sq)=C_{0s}+A_{0s}\in\Real\oplus\bar\h.$$
This shows that $\bar C_{0s}+A_{0s}$ preserves $L$. It holds
$A_{0s}=I_sT_0+T_0I_s$. Since $T_0$ takes values in $L$, for any
$X\in L$ it holds
\begin{equation}\label{pri1}\pr_{\spa_\H L'}(\bar C_{0s}X+I_sT_0X)\in
L'.\end{equation} Suppose that $L'$ is of type $B(l)$, $l\geq 2$
and it is given by vectors $f_1,...,f_l$. Then
$$\pr_{\spa_\H L'}T_0f_1=\pr_{L'}T_0f_1=a_1f_1+\cdots+a_lf_l+b_1(if_1+jf_2)+\cdots+b_{l-1}(if_{l-1}+jf_{l}),$$
where $a_1,...,a_l,b_1,...,b_{l-1}\in\Real$. From \eqref{pri1} it
follows that $\bar C_{0s}f_1+I_s\pr_{L'}T_0f_1\in L'$. Taking
$s=1$, we get
$$\bar C_{01}f_1+a_1if_1+\cdots+a_lif_l-b_1f_1\cdots-b_{l-1}f_{l-1}+kb_1f_2+\cdots+kb_{l-1}f_l\in
L'.$$ Hence, $b_1=\cdots=b_{l-1}=a_2=\cdots =a_l=0$ and $\bar
C_{01}=c_1-a_1i$ for some $c_1\in\Real$. In particular,
$\pr_{L'}T_0f_1=a_1f_1$. Similarly, we get $\bar C_{02}=c_2-a_1j$
and $\bar C_{03}=c_3-a_1k$ for some $c_2,c_3\in\Real$. Using the
equality $C_{03}=C_{02}i-C_{01}j$, we obtain $a_1=c_1=c_2=c_3=0$.
Hence, $C_{rs}=0$. This shows that $\bar\h\subset\sp(n)$ and
$\pr_\H\g=0$. For $L'$  of type $A(2l-1)$, $l\geq 2$, the proof is
similar, hence $\bar\h\subset\sp(n)$ for any $L'$. Let
$\h=\bar\h$. We claim that $\h$ preserves decomposition
\eqref{decLvdokve}. Since $\h$ commutes with $I_1,I_2,I_3$, it
preserves $\H^m=L\cap I_1 L\cap I_2 L\cap I_3 L$. Hence $\h$
preserves $(\H^m)^{\bot_\eta}=L_2=\im
\H^{m_1}\oplus\Co^{m_2}\oplus L'$. Next, $\h$ preserves $L_2\cap j
L_2=i\Real^{m_1}\oplus k\Real^{m_1}$ and $L_2\cap k
L_2=i\Real^{m_1}\oplus j\Real^{m_1}$, i.e. it preserves
$i\Real^{m_1}$. Thus $\h$ preserves $\Real^{m_1}$ and $\im
\H^{m_1}$. By similar arguments, $\h$ preserves $\Co^{m_2}$ and
$L'$. The claim is proved.

 The space
$\H^n=\spa_\H L$ is the direct sum of four subspaces and $\h$
preserves this decomposition, hence we may write
$A_{rs}=A^1_{rs}+A^2_{rs}+A^3_{rs}+A^4_{rs}$. Similarly, we
decompose the elements $S_{rs}$ and $T_s$. Equality \eqref{p1To}
shows that $T^1_s=T_s|_{\H^m}$, $T^2_s=T_s|_{\H^{m_1}}$,
$T^3_s=T_s|_{\H^{m_2}}$ and $T^4_s=T_s|_{\spa_\H L'}$. Clearly,
these maps take values in $\H^m$, $\im \H^{m_1}$, $\Co^{m_2}$ and
$L'$, respectively. Since $\h$ preserves each summand in the
direct sums $\Co^{m_1}\oplus j\Co^{m_1}$, $\Co^{m_2}\oplus
j\Co^{m_2}$ and $L'\oplus iL'$, and acts in each summand
simultaneously, according to \cite{ESI}, $R'\in\R(\h\cap\sp(m))$
and $P_0\in\R(\h\cap\sp(m))$.

Let $S^2_{01}=is_1+js_2+ks_3$, $S^2_{02}=is_6+js_4+ks_5$, where
$s_1,...,s_6\in\Real^{m_1}$. The condition
$S^2_{03}=jS^2_{01}-iS^2_{02}\in\Im\H^{m_1}$ is equivalent to the
equality $s_6=s_2$. The vectors $s_1,...,s_5\in\Real^{m_1}$ are
arbitrary.

Since $S^3_{01},S^3_{02},S^3_{03}\in\Co^{m_2}$, and
$S^3_{03}=jS^3_{01}-iS^3_{02}$, we see that $S^3_{01}=0$ and
$S^3_{02}\in\Co^{m_2}$ may be arbitrary.

Since $A_{rs}$ preserves $\Im\H^{m_1}$ and $A_{rs}\in\sp(n)$, it
preserves $i\Im\H^{m_1}\cap j\Im\H^{m_1}\cap
k\Im\H^{m_1}=\Real^{m_1}$. Let $X\in\Real^{m_1}$, then
$$T^2_1(X)=-\frac{1}{2}(A^2_{01}(X)+I_3A^2_{02}(X)-I_2A^2_{03}(X))\in\Im\H^{m_1},$$
hence, $A^2_{01}(X)=0.$ Since $A^2_{01}\in\sp(m_1)$, this implies
$A^2_{01}=0$. Similarly, $A^2_{02}=A^2_{03}=0$.

Let $X\in\Co^{m_2}$. Since
$$T^3_0(X)=-\frac{1}{2}(I_1A^3_{01}(X)+I_2A^3_{02}(X)+I_3A^3_{03}(X))\in\Co^{m_2},$$
$I_2A^3_{02}(X)+I_3A^3_{03}(X)=I_2(A^3_{02}(X)-I_1A^3_{03}(X))$
and $A^3_{0r}(X)\in \Co^{m_2}$ for any $r$,  we get that
$A^3_{02}(X)=I_1A^3_{03}(X)$. This implies
$A^3_{02}|_{\Co^{m_2}}=I_1A^3_{03}|_{\Co^{m_2}}$. Hence,
$A^3_{02}|_{\Co^{m_2}}=A^3_{03}|_{\Co^{m_2}}=0$  and
$A^3_{02}=A^3_{03}=0$. Next,
$T^3_2(X)=I_2T^3_0(X)-A^3_{02}(X)=\frac{1}{2}I_3 A^3_{01}(X)$.
Hence, $A^3_{01}(X)=0$. This shows that $A^3_{rs}=0$.

Let $Y\in  L'$. Then for $s\neq 0$,
$T^4_s(Y)=I_sT^4_0(Y)-A^4_{0s}(Y)\in L'$. Since $L'\cap I_s L'=0$
and $A^4_{0s}$ preserves $L'$, we get $T^4_0(Y)=0$. From
\eqref{(*)} applied to the vectors $I_s q$, $X,Y\in\H^n$, $q$,  it
follows that
$$\eta(T_s(X),Y)=\eta(I_sT_0(Y),X).$$
Let $Y\in L'$, we get $\eta(T^4_s(X),Y)=0$ for any $X\in\H^n$ and
any $s$. Hence $T^4_s$=0. Consequently, $A^4_{rs}=0$. Thus,
$\h\subset\sp(m)$.

We see that $L'$ must be spanned by  elements
$S_{01},S_{02},S_{03}\in L'$  that satisfy
$S_{03}=jS_{01}-iS_{02}$, i.e. $L'$ must satisfy
\begin{equation}\label{cond3} L'=\rho(L'),\quad\text{where}\quad
\rho(L')=\spa_\Real\{X,Y,jX-iY|\, X,Y,jX-iY\in L'\}.\end{equation}
Clearly, the space $B(l)$, $l\geq 2$ satisfies this condition,
while the space $B(1)$ does not satisfy this condition.

\begin{lemma} The space $L'=A(2l-1)$, $l\geq 2$ does not satisfy
the condition \eqref{cond3}. \end{lemma}

{\it Proof.} It can be directly checked that $\rho(A(3))=0$.
 We claim that if $L'=A(2l-1)$, $l\geq 3$,
then
\begin{multline*}\rho(L')=\spa_\Real\{f_1,...,f_{l-1},f_{l+1},...,f_{2l-1},if_1+jf_2,...,if_{l-2}+jf_{l-1},
jf_{l+1}+if_{l+2},...,jf_{2l-2}+if_{2l-1}\}.\end{multline*}  We
prove this claim using the induction over $l$. For $l=3$ this can
be checked directly. Suppose that the claim holds for some $l\geq
3$. We will prove it for $l+1$.

Clearly, $A(2(l+1)-1)$ can be obtained from $A(2l-1)$ adding some
vectors $f_0,f_{2l}$, $if_0+jf_1$, $jf_{2l-1}+if_{2l}$. Let
$X,Y\in A(2(l+1)-1)$. Then,
\begin{align*}X&=af_0+bf_{2l}+c(if_0+jf_1)+d(jf_{2l-1}+if_{2l})+\tilde
X,\\ Y&=xf_0+yf_{2l}+u(if_0+jf_1)+v(jf_{2l-1}+if_{2l})+\tilde Y
\end{align*} for some $a,b,c,d,x,y,u,v\in\Real$, $\tilde X,\tilde
Y\in A(2l-1)$. It can be checked that if $jX-iY\in A(2(l+1)-1)$,
then $a,b,c,d=0$ and  $$j X-i Y=j\tilde X_1-i\tilde
Y-x(if_0+jf_1)-y(jf_{2l-1}+if_{2l})-uf_2-vf_{2l-2}+uf_0+vf_{2l},$$
 where
$$\tilde X_1=\tilde
X+xf_1-yf_{2l-1}-u(if_1+jf_2)-v(jf_{2l-2}+if_{2l-1})\in A(2l-1).$$
Hence, $j\tilde X_1-i\tilde Y\in A(2l-1)$. This and the induction
hypothesis prove the inclusion $\subset$, the inverse inclusion is
obvious. The lemma is proved. \hfill $\square$

Thus $L'$ is an $g$-orthogonal sum of the spaces of the form
$B(l)$, $l\geq 2$.

{\bf 2.a)} Suppose that $L'=0$, $m_1\neq 0$ and $m_2= 0$, i.e.
$L=\H^m\oplus\im \H^{m_1}$, $m+m_1=n$. For simplicity of the
exposition we may assume that $m=0$, i.e. $L=\im \H^n$. Obviously,
$\tilde\h\cap\sp(1)=0$ and elements of the form $a-\Op(aE_n)$
(where $a\in\sp(1)$) preserve $L=\im \H^n$, hence
$\tilde\h\subset\{a-\Op(aE_n)|\,a\in\sp(1)\}\oplus\h_1$, where
$\h_1\subset\sp(n)$ is a vector subspace preserving $\im \H^n$.
Let $R\in\R(\g)$ be as in Proposition \ref{propR}. Then
$$\pr_{\H\oplus\sp(n)}R(q,I_sq)=C_{0s}+A_{0s}=C_{0s}+\Op(a_{0s}E_n)+B_{0s},$$
where $a_{0s}=\im C_{0s}$ and $B_{0s}\in\sp(n)$ preserves $\im
\H^n$. Clearly, $B_{0s}$ preserves $\Real^n$. Recall that
$$T_0(X)=-\frac{1}{2}(I_1A_{01}X+I_2A_{02}X+I_3A_{03}X)\in L$$
for any $X\in\H^n$. Let $e_\alpha\in\Real^n\subset\H^n$ be an
element of the basis. The condition $T_0(e_\alpha)\in L$ implies
$\Re(ia_{01}+ja_{02}+ka_{03})=0$. The condition $T_0(ie_\alpha)\in
L$ implies
$$B_{01}e_\alpha=\Re((ia_{01}+ja_{02}+ka_{03})i)e_\alpha.$$
Since $B_{01}\in\sp(n)\subset\so(4n)$, we conclude
$\Re((ia_{01}+ja_{02}+ka_{03})i)=0$ and $B_{01}=0$. Similarly,
$B_{02}=B_{03}=0$ and
$$\Re((ia_{01}+ja_{02}+ka_{03})j)=\Re((ia_{01}+ja_{02}+ka_{03})k)=0.$$
Thus, $ia_{01}+ja_{02}+ka_{03}=0$, i.e. $a_{03}=ja_{01}-ia_{02}$.
This and the equalities $a_{0s}=\im C_{0s}$,
$C_{03}=C_{02}i-C_{01}j$ imply $a_{0s}=C_{0s}\in\im \H$. Thus
since $\g$ is a Berger algebra, $\pr_\Real\g=0$. Lemma \ref{L1}
shows that either $\bar\h=0$, or
$\bar\h=\{a+\Op(aE_{n})|\,a\in\sp(1)\}$. If we do not assume that
$m=0$, then   $\bar\h\subset
\{a+\Op(aE_{m_1})|\,a\in\sp(1)\}\oplus\sp(m)$ and $\pr_\Real\g=0$,
moreover, the projection of $\bar\h$ to
$\{a+\Op(aE_{m_1})|\,a\in\sp(1)\}$ is either trivial or it
coincides with $\{a+\Op(aE_{m_1})|\,a\in\sp(1)\}$.

{\bf 2.b)} Suppose that  $m_1=0$ and $m_2\neq 0$, i.e.
$L=\H^m\oplus\Co^{m_2}$, $m+m_2=n$. As above, suppose that $m=0$,
then $L=\Co^n$.  Let $e_1,...,e_n$ a basis of the complex vector
space $L=\Co^n$. Let $C+A\in\bar\h$, where $C\in\sp(1)$ and
$A\in\sp(n)$. Let $A_{\alpha\beta}$ be the matrix of $A$ with
respect to the basis $e_1,...,e_n$ of $\H^n$. Then since
$(-C+A)e_\alpha\in L$ and $(-C+A)ie_\alpha\in L$, we get
$-C+A_{\alpha\alpha}\in\Co$ and $-Ci+iA_{\alpha\alpha}\in\Co$.
Consequently, $C\in\Co$. This shows that $A$ preserves $L$ and
$\bar\h\subset\Real i\oplus\sp(n)$. Let $R\in\R(\g)$ be as in
Proposition \ref{propR}. Then
$$\pr_{\H\oplus\sp(n)}R(q,I_sq)=C_{0s}+A_{0s}\in\Co\oplus\sp(n)$$
 preserves $L$ and $A_{0s}$ preserves $L$. By the arguments of the
 proof of
 statement 1), $A_{0s}=0$. Thus, $\bar \h\subset \Real i$.
 If $m\neq 0$, then  $\bar \h\subset \Real
i\oplus\sp(m)$.

{\bf 2.c)} Suppose that $m_1\neq 0$ and $m_2\neq 0$, i.e.
$L=\H^m\oplus\im \H^{m_1}\oplus\Co^{m_2}$, $m+m_1+m_2=n$. As in
the proof of 2.b), it can be shown that $\bar\h\subset\Real
i\oplus\sp(m)\oplus\sp(m_1)$, i.e. $\pr_{\sp(1)}\g\subset i\Real$.
As in the proof of 2.a), it can be proved that $\pr_{\Real}\g=0$.
From this and Lemma \ref{L1} it follows that $\pr_{\H}\g=0$, i.e.
$\bar\h\subset\sp(m)\oplus\sp(m_1)$. By the arguments of the proof
of 1), $\bar\h\subset\sp(m)$.  The lemma is proved. \hfill
$\square$

Now using Theorem \ref{ThRes31} and Lemmas \ref{L1}, \ref{LL'}, it
is easy to obtain the classification of weakly irreducible Berger
subalgebras of $\sp(1,n+1)_{\H p}$ containing $\Im\H$. All such
subalgebras are exhausted by the Lie algebras given in the
statement of the Theorem \ref{Clas}. Let us consider some
examples.

Let $L=\H^{m}\oplus\Im\H^{m_1}\oplus\Co^{m_2}\oplus L'$ and
$L'\neq 0$. In this case, for any weakly irreducible Berger
subalgebra $\g\subset\sp(1,n+1)_{\mathbb{H}{p}}$  by Lemma
\ref{LL'}, we have $\bar\h\subset\sp(m)$ and $\pr_\H\g=0$. Then
$\g$ can be of Type $\varphi$ or Type $\psi$ from Theorem
\ref{ThRes31}. If $\g$ is of Type $\varphi$, then we get that
$\g=\g_{II}$ from Theorem \ref{Clas} with $\varphi=0$. If $\g$ is
of Type $\psi$, then $\g=\g_{IV}$. Next we may assume that $L'=0$.

Let $L=\H^{m}\oplus\Im\H^{m_1}$. By Lemma \ref{LL'},
$\pr_\Real\g=0$ and either $\bar{\h}= \{a+\Op(aE_{m_1})|\,a\in\im
\H\}\oplus\h$, or $\bar{\h}=\{\phi(A)+\Op(\phi(A)E_{m_1})+A|\,
A\in\h\},$ where $\h\subset\sp(m)$ is a subalgebra and
$\varphi:\h\to\sp(1)$ is a homomorphism. Using the fact that there
are no two-dimensional subalgebras of $\sp(1)$ and Lemma \ref{L1},
we get that either $\phi=0$, or $Im\phi=\sp(1)$. Since
$\pr_\Real\g=0$,  $\g$ is either of Type $\varphi$ (with
$\varphi=0$) or of Type $\psi$. In the first case $\g$ coincides
with $\g_{III1}$ or $\g_{III2}$.

Let $\g$ be of Type $\psi$. Since $\tilde h$ annihilates the
subspace $U\subset L$, it can not contain $\sp(1)$, i.e.
$\bar{\h}=\{\phi(A)+\Op(\phi(A)E_{m_1})+A|A\in\h\}.$ If $\phi=0$,
then  $\g=\g_{IV}$.

Assume that $\phi\neq 0$, then $\phi$  is surjective. Now we find
subspaces $U\subset L$ such that $\tilde{\h}U=0$. Any $u\in U$ has
the  form $u=u_1+u_2$, where $u_1\in\H^{m}$, $u_2\in\Im\H^{m_1}$.
Let $\xi=-\phi(A)+\Op(\phi(A)E_{m_1})+A\in\tilde{\h}$, then $\xi
u_1=\xi u_2=0$. We see that if $u_2=\sum_{i=m+1}^{n}a_ie_i$, where
$e_{m+1},...,e_{n}$ is a basis of $\Real^{m_1}$, $a_i\in\Im\H$,
then $\xi u_2=-\sum_{i=m+1}^{n}[\phi(A),a_i]e_i=0$. From here
$a_i=0$ for all $i=m+1,...,n$, i.e. $u_2=0$. Thus, $u\in\H^{m}$
and $U\subset\H^m$.

We claim that for any $\xi\in\sp(1)$, $\xi\neq 0$ it holds $\xi
U\cap U=0$. Indeed, let $\xi U\cap U\neq 0$, then there exists a
non-zero vector $u\in U$ such that $\xi u\in U$. Since
$\tilde{\h}$ annihilates $u$ and $\xi u$,  $$(\phi
(A)-\Op(\phi(A)E_{m_1})-A)u=0,\quad (\phi
(A)-\Op(\phi(A)E_{m_1})-A)\xi u=0$$ for all $A\in\h$. This implies
that $[\xi,\phi(A)]u=0$ hence, $[\xi,\phi (A)]=0$ for any
$A\in\h$.
 Consequently, $[\xi,\sp(1)]=0$, where
$\xi\in\sp(1)$. It gives the contradiction. Thus, according to
\cite{DR78}, $U$ is a $g$-orthogonal direct sum of vector spaces
of the form $A(2l-1)$ and $B(l)$, $l\geq 1$. Further, assume that
$U$ contains a subspace of the form $A(2l-1)$ or $B(l)$, where
$l\geq 2$. Then $f_1, f_2, if_1+jf_2\in U$ and
\begin{multline*}(\phi (A)-\Op(\phi(A)E_{m_1})-A)f_1=0, \quad (\phi
(A)-\Op(\phi(A)E_{m_1})-A)f_2=0, \\(\phi
(A)-\Op(\phi(A)E_{m_1})-A)(if_1+jf_2)=0\end{multline*} for all
$A\in\h$. Hence, $Af_1=\phi (A)f_1$ and $Af_2=\phi (A)f_2$. We get
$\phi (A)i=i\phi (A)$ and $\phi (A)j=j\phi (A)$ for any $A\in\h$.
Consequently, $\phi=0$. This gives a contradiction. Thus $U$ is a
$g$-orthogonal sum of vector subspaces of the form $B(1)$.
Consider $\spa_{\H} U\subset\H^{m}$ and denote  by $\H^{m_0}$ its
orthogonal complement. Then $\H^{m}=\H^{m_0}\oplus\spa_{\H} U$.
Let $e_1,...,e_{m_0}$ be a $g$-orthonormal basis of $\H^{m_0}$.
Since $U$ is a direct sum of subspaces of the form
$B(1)=\Real{f_1}$,
$U=\spa_{\Real}\{e_{m_{0}+1},...,e_{m}\}=\Real^{m-m_0}$, and the
vectors $e_{m_{0}+1},...,e_{m}$ are $g$-orthogonal.

Now let us find the algebras $\tilde{\h}$ that annihilate $U$. Let
$\phi (A)-\Op(\phi(A)E_{m_1})-A\in\tilde \h$. Then $(\phi
(A)-A)e_i=0$, $i=m_{0}+1,...,m$, i.e. $Ae_i=\phi (A)e_i$. This
shows that $A|_{\H^{m-m_0}}=\Op(\phi(A)E_{m-m_0})$. Moreover, $A$
preserves $\H^{m-m_0}$, and, consequently, it preserves
$\H^{m_0}$. It is clear that the obtained properties of $\tilde\h$
are equivalent to the conditions
$$\tilde\h\subset \{a-\Op(aE_{m-m_0})-\Op(aE_{m_1})|\, a\in\sp(1)\}\oplus\sp(m_0),
\quad \tilde h\not\subset\sp(m_0),$$ or to the conditions
$$\bar\h\subset \{a+\Op(aE_{m-m_0})+\Op(aE_{m_1})|\,a\in\sp(1)\}\oplus\sp(m_0),
\quad \bar h\not\subset\sp(m_0).$$ Suppose that $\bar\h=
\{a+\Op(aE_{m-m_0})+\Op(aE_{m_1})|\, a\in\sp(1)\}\oplus\h$, where
$\h\subset\sp(m_0)$. Since $\psi|_{\bar h'}=0$, $\psi$ is zero on
the first summand. We see that $\g=\g_{V1}$. If $\bar\h$ does not
contain $\{a+\Op(aE_{m-m_0})+\Op(aE_{m_1})|\, a\in\sp(1)\}$, then
$\g=\g_{V2}$.

Further, let $L=\H^{m}\oplus\Co^{m_2}$. By Lemma \ref{LL'},
$\bar{\h}\subset\Real i\oplus\sp(m)$. Assume that $\bar{\h}=\Real
i\oplus\h$, $\h\subset\sp(m)$. In this case, $\g$ can be of Type
$\Real$ or Type $\varphi$. Let $\g$ be of Type $\Real$, then
$\g=\g_{I1}$. If $\g$ is of Type $\varphi$ and $\varphi\neq 0$, we
have $\{\varphi (A)+A |\, A\in\bar{\h}\}=\Real (\varphi
(i)+i)\oplus\{\varphi (A)+A |\, A\in\h\}$, $\h\subset\sp(m)$. If
$\varphi (i)=0$, then $\g=\g_{I3}$. If $\varphi (i)\neq 0$, then
$\g=\g_{I5}$. Suppose that $\bar{\h}=\{\varphi (A)i+A |\,
A\in\h\}$, $\h\subset\sp(m)$. Again, $\g$ can be of Type $\Real$
or Type $\varphi$. In the first case, $\g=\g_{I2}$. In the second
case
$$\{\varphi (A)+A |\, A\in\bar{\h}\}=\{\varphi (\phi (A)i+A)+\phi
(A)i+A |\, A\in\h\}=\{\widetilde{\varphi} (A)+\phi (A)i+A |\,
A\in\h\},$$ here $\widetilde{\varphi} (A)=\varphi (\phi (A)i+A)$,
$\h\subset\sp(m)$.  By Lemma \ref{L1}, $\widetilde{\varphi}$ and
$\phi$ are not proportional. Thus, $\g=\g_{I4}$. Now let
$\bar{\h}\subset\sp(m)$. Then $\g$ can be of Type $\varphi$. Using
Lemma \ref{L1}, we get that $\g=\g_{II}$ with $\varphi=0$. Also,
$\g$ can be of Type $\psi$, then $\g=\g_{IV}$.

All other cases can be considered in the similar way.

{\bf Any weakly irreducible Berger subalgebra
$\g\subset\sp(1,n+1)_{\H p}$ contains $\im\H$}

\begin{prop} Let $\g\subset\sp(1,n+1)_{\H p}$  be
a weakly irreducible  Berger subalgebra, then it is conjugated to
a subalgebra that contains the ideal $\Im\H$.
\end{prop}

{\bf Proof.} Let $\g\subset\sp(1,n+1)_{\H p}$ be a weakly
irreducible subalgebra. Let $f:\sp(1,n+1)_{\H p}\to\simil\H^n$ be
the homomorphism as in Section \ref{w}. According to \cite{B1} and
Section \ref{w}, the image $f(\g)\subset\simil\H^n$ coincides with
one of the following algebras:
\begin{description}
\item{{\bf Type $\mathbb{R}$.}} $f(\g)=(\mathbb{R}\oplus \tilde\h)\ltimes
L$, \item{{\bf Type $\varphi$.}}
$f(\g)=\{\varphi(A)+A|\,\,A\in\tilde\h\}\ltimes L$,
\item{{\bf Type $\psi$.}} $f(\g)=\{A+\psi(A)|\,\,A\in\tilde\h\}\ltimes W$,
\end{description}
where $L$ is as in \eqref{decLvdokve}, $\tilde\h\subset
\sp(1)\oplus\sp(n)$ is a subalgebra, $\tilde\h L\subset L$;
$\varphi\in\Hom(\tilde\h,\mathbb{R})$, $\varphi|_{\tilde\h'}=0$;
for the last algebra we have an orthogonal decomposition
$L=W\oplus U$, $\tilde\h$ annihilates $U$ and $\psi:\tilde\h\to U$
is surjective linear map with $\psi|_{\tilde\h'}=0$.

Suppose that $\g$ is a Berger algebra. The structure of the Lie
brackets of $\sp(1,n+1)_{\H p}$ (see Section \ref{REZ}) shows that
if $m\neq 0$, or $m_1\neq 0$, then $\g$ contains $\Im\H$. Thus,
$L=\Co^{m_2}\oplus L'$. By Lemma \ref{LL'}, if $\g$ is a Berger
algebra with such $L$, then $\pr_{\sp(n)}\g=0$. Hence, if $L'=0$,
then $\bar\h\subset\Real i$; if $L'\neq 0$, then $\bar\h=0$.
Recall that since $\g$ is a Berger algebra, $\dim\pr_{\H}\g\neq
1$. This shows that either $f(\g)$ is of Type $\mathbb{R}$ with
$\bar\h=\Real i$ and $L'=0$, or $f(\g)$ is of Type $\varphi$ with
$\varphi=0$ and $\bar\h=0$.

Consider the first case. We have $(1,0,0,b)\in\g$ for some
$b\in\Im\H$. Let $X\in\Co^n$, then $(0,0,X,c)\in\g$. Next,
$$[(1,0,0,b),(0,0,X,c)]=(0,0,X,2c)\in\g.$$
This shows that $\Co^n\subset\g$. Consequently,
$$[(0,0,e_1,0),(0,0,ie_1,0)]=(0,0,0,-2i)\in\g.$$ Hence,
$\Real(0,0,0,i)\subset\g$. If $(0,0,0,\alpha j+\beta k)\in\g$ for
some $\alpha,\beta\in\Real$ with $\alpha^2+\beta^2\neq 0$, then
taking the Lie bracket of this element with $(i,0,0,c)\in\g$, we
get $(0,0,0,\alpha k-\beta j)\in\g$, i.e. $\Im\H\subset\g$. Assume
that $\g\cap\Real j\oplus\Real k=0$. Then it is not hard to see
that $$\g=\Real(1,0,0,\alpha j+\beta k)\oplus\Real(i,0,0,-\beta
j+\alpha k)\oplus\Co^n\zr\Real(0,0,0,i).$$ Let $R\in\R(\g)$. Above
we have seen that the elements defining $R$ are zero possibly
except for $C_{02}\in\Co$, $S_{02}\in \Co^n$ and some of $D_{rs}$.
Let $X\in\H^n$. It holds $R(I_sq,X)=(0,0,0,\theta_s(X))\in\g$.
Hence, $\theta_s(X)\in\Real i$. From \eqref{LP} it follows that
$\eta(\theta_s(X)p,I_2q)=\eta(I_sS_{02},X)$. Consequently,
$S_{02}=0$, and $\g$ is not a Berger algebra.

Suppose that $f(\g)$ is of Type $\varphi$ and $m_2\neq 0$. Then
for some $b,c\in\Im\H$, $(0,0,e_1,b)$, $(0,0,ie_1,c)\in\g$. Taking
the Lie brackets of these elements, we get $(0,0,0,i)\in\g$. Using
\eqref{LP}, we obtain that
$\eta(\theta_s(X)p,I_2q)=\eta(I_sS_{02},X)$. If $\g$ is a Berger
algebra, then for some $R\in\R(\g)$ it holds $S_{02}\neq 0$. There
exists $X\in\Im\H$ such that $\eta(S_{02},X)\neq 0$. Hence,
$\theta_0(X)\in\Im\H$ has a non-zero projection to $\Real
j\subset\Im\H$. Since $R(I_sq,X)=(0,0,0,\theta_s(X))\in\g$, there
exists $\alpha,\beta\in\Real$ such that $\alpha\neq 0$ and
$(0,0,0,\alpha j+\beta k)\in\g$. We may assume that
$\alpha^2+\beta^2=1$. In \cite[Lemma 6]{44} it is shown that there
exists $x,y\in\Real$ such that $x^2+y^2=1$ and with respect to the
new basis with $p'=(x+iy)p$ and $q'=(x+iy)q$ the elements
$(0,0,0,i)\in\g$ and $(0,0,0,\alpha j+\beta k)\in\g$ have the form
$(0,0,0,i)$ and $(0,0,0,j)$, respectively. Note that
$S_{03}=-iS_{02}\neq 0$. As above,  there exists $X\in\Im\H$ such
that $\eta(I_1S_{03},X)\neq 0$. Hence,
$\eta(\theta_s(X)p,I_3q)\neq 0$, i.e. $\theta_1(X)\in\Im\H$ has a
non-zero projection to $\Real k\subset\Im\H$. We conclude that
$\g$ contains $\Im\H$.

Suppose now that $m_2=0$, i.e. $L=L'$. Suppose that $\dim
\g\cap\Im\H=2$. From the \cite[Lemmas 5, 6]{44} it follows that
choosing an appropriate basis we may get $\g\cap\Im\H=\Real
i\oplus\Real j$. Since $R(I_sq,X)=(0,0,0,\theta_s(X))$, we see
that $\theta_s(X)\in \Real i\oplus\Real j$ for any $X\in\H^n$.
Since $\theta_s(X)=g(X,S_{0s})+I_s\theta_0(X)$, we obtain that for
any $X\in\H^n$ and $s=1,2$ it holds
$g(X,S_{0s})\in\Real\oplus\Real i\oplus\Real j$. Since
$g(kS_{0s},S_{0s})=kg(S_{0s},S_{0s})\in\Real k$, we have
$g(S_{0s},S_{0s})=0$, consequently, $S_{0s}=0$ for $s=1,2$. This
implies $S_{rs}=0$. Hence, $\g$ is not a Berger algebra. The case
$\dim \g\cap\Im\H<2$ follows from this one. This proves the
proposition and the theorem.  \hfill $\square$ $\square$


\section{Pseudo-hyper-K\"ahlerian symmetric spaces of index 4}\label{secSym}
 In \cite{AC01,KO07,KO08} indecomposable simply connected
pseudo-hyper-K\"ahlerian symmetric spaces of signature $(4,4n+4)$
are classified. Here we use the results of this paper to give new
proof to this result. For $n=0$ such new proof is obtained in
\cite{44}.


As it is explained e.g. in \cite{AlIRMA} the classification of
indecomposable simply connected pseudo-hyper-K\"ahlerian symmetric
spaces is equivalent to the classification of pairs $(\g,R)$ ({\it
symmetric pairs}), where $\g\subset\sp(1,n+1)$ is a subalgebra,
$R\in\R(\g)$, $R(\Real^{4,4n+4},\Real^{4,4n+4})=\g$, and for any
$\xi\in\g$ it holds
\begin{equation}\label{AR=0}\xi\cdot R=0,\quad (\xi\cdot R)(x,y)= [\xi,R(x,y)]-R(\xi x,y)
-R(x, \xi y),\end{equation} where  $x,y\in \Real^{4,4n+4}$. {\it
An isomorphism }of  symmetric pairs $f:(\g_1,R_1)\to (\g_2,R_2)$
consists of an isometry of $\Real^{4,4n+4}$ that defines the
equivalence of the representations $\g_1,\g_2\subset\sp(1,n+1)$
and sends $R_1$ to $R_2$.
 For a positive real number
$c\in\Real$, the symmetric pairs $(\g,cR)$ and $(\g,R)$ define
diffeomorphic simply connected symmetric spaces and the metrics of
these spaces differ by the factor $c$. Hence we may identify
$(\g,cR)$ and $(\g,R)$.

\begin{theorem} Let $(M,h)$ be a non-flat simply connected pseudo-hyper-K\"ahlerian symmetric
space of signature $(4,4n+4)$ ($n\geq 1$) and
$\g\subset\sp(1,n+1)$ its holonomy algebra. Then $n=2$,
$$\g=L'\zr\Im\H,$$ there exists a basis $e_1,e_2$ of $\H^2$ such
that $L'=\spa_\Real\{e_1,e_2,je_2+ie_1\}$, the Gram matrix of $g$
with respect to this basis equals to
$G=\left(\begin{array}{cc}1&-\frac{1}{2}k\\\frac{1}{2}k&1\end{array}\right).$
The manifold $(M,h)$ is defined by the symmetric pair $(\g,R)$,
where $R$ is defined as in Proposition \ref{propR} and it is given
by $S_{01}=e_1$, $S_{02}=-e_2$ with other elements defining $R$
being zero.
 \end{theorem}

{\bf Proof.} Since $(M,h)$ is Ricci-flat, its holonomy algebra
$\g$ cannot be reductive \cite{AlIRMA}, hence it is conjugated to
a subalgebra of $\sp(1,n+1)_{\H p}$. Let $(\g,R)$ be a symmetric
pair. Then the tensor $R$  is given as in Proposition~\ref{propR}.

First suppose that $\pr_{\H}\g\neq 0$. Let $\xi=(a,A,0,0)\in\g$.
Then $$[\xi,R(I_rq,I_sq)]-R(\xi I_rq,I_sq) -R(I_rq, \xi I_sq
)=0.$$ Taking the projection on $\H$, we get the same equations on
$a$, $C_{rs}$ as in \cite[Section 4]{44}, where it is shown that
these equations imply $C_{rs}=0$. Hence,
$R(\Real^{4,4n+4},\Real^{4,4n+4})\neq\g$ and we obtain a
contradiction. Thus, $\pr_{\H}\g=0$.

Since $\pr_{\H}\g=0$, we get $C_{rs}=B_{rs}=0$. Let
$\xi=(0,0,Y,0)\in\g$ and $X,Z\in\H^n$. The condition $(\xi\cdot
R)(X,Z)=0$ implies $R'=0$ and $\Im g(Y,Q(X,Z))=0$.  This shows
that $Q=0$ and $P_0=0$. Suppose that $Y\in\H^m$. The condition
$(\xi\cdot R)(q,I_sq)=0$ implies
$$-A_{0s}Y+T_s(Y)-T_{0}(I_sY)=0.$$ Substituting
$T_s(Y)=I_sT_0(Y)-A_{0s}(Y)$, we get
$$-2A_{0s}Y+I_sT_0(Y)-T_0(I_sY)=0.$$ Replacing $Y$ by $I_sY$,
multiplying the obtained equality by $I_s$, and combining it with
the last one, we obtain $A_{0s}Y=0$. Hence, $A_{0s}=0$.

Now $R$ is defined only by $S_{rs}$ and $D_{rs}$ and $L$ must be
spanned by $S_{01},S_{02},S_{03}$. This shows that $\g=L\zr\Im\H$
and $L$ has dimension at most 3. Hence $n=1$ or 2. If $n=1$, then
either $L=\Im\H$, or $L=\Co$. If $n=2$, then $L=L'$ and $\dim
L'=3$.

Let $Y\in L$ and $\xi=(0,0,Y,0)\in\g$. The condition $(\xi\cdot
R)(q,I_sq)=0$ implies \begin{equation}\label{esymS}2\Im
g(Y,S_{0s})=\theta_0(I_sY)-\theta_s(Y).\end{equation} From this
equation, \eqref{p1theta0} and \eqref{p1thetas}, we have $$\Im
g(Y,S_{0s})=\frac{1}{2}(I_1g(I_sY,S_{01})+I_2g(I_sY,S_{02})+I_3g(I_sY,S_{03})).$$
Using \eqref{p1S}, we get \begin{equation}\label{AA} \Im
g(Y,S_{0s})=\frac{1}{2}(I_1g(I_sY,S_{01})+I_2g(I_sY,S_{02})-
I_3g(I_sY,S_{01})I_2+I_3g(I_sY,S_{02})I_1),\end{equation} where
$Y\in L$. Substituting in turn $Y=S_{01}$ and $s=1$; $Y=S_{02}$
and $s=2$ in the last equation, we obtain the following
\begin{align}\label{AA1}0&=-kg(S_{01},S_{02})+jg(S_{01},S_{02})i, \\
\label{AA2} 0&=kg(S_{02},S_{01})+ig(S_{02},S_{01})j.\end{align}
Taking the conjugation to \eqref{AA2}, we get
\begin{equation}\label{AA3}
0=-g(S_{01},S_{02})k+jg(S_{01},S_{02})i.\end{equation} From
\eqref{AA1} and \eqref{AA3}, we have
$g(S_{01},S_{02})k-kg(S_{01},S_{02})=0$. It means that
$g(S_{01},S_{02})=\alpha+\beta k$, $\alpha,\beta\in\Real$.
Substituting this to \eqref{AA1}, we obtain $\alpha=0$, i.e.
\begin{equation} \label{AA4}g(S_{01},S_{02})=\beta k, \quad
\beta\in\Real.\end{equation} If we take $Y=S_{01}$ and $s=2$ in
\eqref{AA}, then
\begin{equation}\label{AA5}g(S_{01},S_{01})=2\beta.\end{equation}
Similarly, $g(S_{02},S_{02})=2\beta$. If $\beta=0$, then
$S_{01}=S_{02}=0$ and consequently, $S_{rs}=0$.

Now let $\beta\neq 0$, then $S_{01}\neq 0$ and $S_{02}\neq 0$. Let
$n=1$. If $L=\Co e$, then $S_{01}, S_{02}, S_{03}\in\Co e$. At the
same time $S_{03}=jS_{01}-iS_{02}\in\Co e$, from here $S_{01}=0$
and $\beta=0$. Further, let $L=\Im\H e$. Since $S_{01}\neq 0$,
there exists $a\in\H$ such that $S_{02}=aS_{01}$. Using
\eqref{AA4} and \eqref{AA5}, we get $a=-\frac{1}{2}k$, i.e.
$S_{02}=-\frac{1}{2}kS_{01}$ and $S_{03}=\frac{1}{2}jS_{01}$. The
condition $S_{0s}\in\Im\H e$ implies $S_{01}=\gamma ie$,
$\gamma\in\Real$ (see \cite[Section 4]{44}). Therefore,
$S_{02}=-\frac{1}{2}\gamma je$ and $S_{03}=-\frac{1}{2}\gamma ke$.
Substituting $Y=je$, $s=3$ and obtained above $S_{01}, S_{02},
S_{03}$ to \eqref{AA}, we prove that $\gamma=0$.

Let $n=2$, then $L=L'$ and we obtain $L'=\{e_1,e_2,j e_1+i e_2\}$
for some basis $(e_1,e_2)$ of $\H^2$. Let us find $S_{rs}$. One
may write $$S_{0s}=a_se_1+b_se_2+c_s(j e_1+i e_2),$$ where
$a_s,b_s,c_s\in\Real$. The condition $S_{03}=jS_{01}-iS_{02}$
implies $S_{01}=a_1e_1$ and $S_{02}=-a_1e_2$. We may assume that
$a_1=1$. From the above reasoning it follows that
$G=\left(\begin{array}{cc}2\beta&-\beta k\\\beta
k&2\beta\end{array}\right)$ for some $\beta\in\Real$, $\beta>0$.
Changing $e_1, e_2$ by
$\frac{\sqrt{2\beta}}{2\beta}e_1,\frac{\sqrt{2\beta}}{2\beta}e_2,$
we get
$G=\left(\begin{array}{cc}1&-\frac{1}{2}k\\\frac{1}{2}k&1\end{array}\right)$.

Note that we still have possibly non-zero elements $D_{rs}$. Let
us consider $p'=p$,  $q'=-\frac{g(X,X)}{2}p+X+q$ for some
$X\in\H^n$. This defines a new basis $p',e'_1,e'_2,q'$ of
$\H^{1,3}$. With respect to this basis $R$ is given by the
elements $S'_{rs}$ and $D'_{rs}$. As we have  shown just now, the
vectors $e'_1$ and $e'_2$ can be chosen to have the same
properties as $e_1$, $e_2$, i.e. we may assume that $S'_{rs}$
remain the same. It holds
$$D'_{rs}=D_{rs}-\theta_s(I_rX)+\theta_r(I_sX)-g(X,S_{rs})+g(S_{rs},X).$$
Choosing $r=0$ in the above equation and using \eqref{p1thetas},
we get
$$D'_{0s}=D_{0s}+2\theta_0(I_sX)+2\Im g(S_{0s},X).$$
In the proof of Theorem 4 (see \cite[Section 7]{44}), it has been
shown that if we consider the new vectors $\tilde p=xp$, $\tilde
q=xq$, for some $x\in\H$, then $D_{01}=\mu i$, $D_{02}=\lambda j$,
where $\mu, \lambda\in\Real$ (again, this will not change
$S_{rs}$). Now take $X=\frac{\mu}{4}ie_1+\frac{\mu}{4}je_2$, then
$$D'_{01}=0, \quad D'_{02}=(\mu+\lambda)j, \quad
D'_{03}=-(\mu+\lambda)k.$$ Again, by similar reasoning as in
\cite[Section 7]{44}, after some transformation $\tilde p=xp$,
$\tilde q=xq$, $x\in\H$, $x\bar{x}=1$, we obtain
$$D'_{01}=(\mu+\lambda)i, \quad D'_{02}=-(\mu+\lambda)j, \quad
D'_{03}=0.$$ Once more, consider $q'=-\frac{g(X,X)}{2}p+X+q$,
where $X=\frac{\mu+\lambda}{4}ie_1+\frac{\mu+\lambda}{4}je_2$,
then $D'_{01}=D'_{02}=0$. Hence, $D'_{rs}=0$. This proves the
theorem. \hfill $\square$

The symmetric spaces from the above theorem and from \cite[Theorem
4]{44} coincide with the ones obtained in \cite{KO07,KO08}.

\end{document}